\documentclass[a4paper,12pt]{article}
\usepackage{amsmath,amsthm,amssymb}
\usepackage[matrix,arrow,curve]{xy}

\newtheorem{Lemma}{Lemma}
\newtheorem{Def}{Definition}
\newtheorem{prop}{Proposition}

\newtheorem{theorem}{Theorem}
\newtheorem{coro}{Corollary}

\newtheorem{rema}{Remark}
\newcommand{\Id}{{\rm Id}}
\newcommand{\charr}{{\rm char\,}}
\newcommand{\ord}{{\rm ord\,}}
\DeclareMathOperator{\Ker}{Ker}
\DeclareMathOperator{\Cok}{Cok}
\DeclareMathOperator{\Aut}{Aut}

\renewcommand{\Im}{{\rm Im\,}}

\newcommand{\Mod}{{\rm Mod}}

\newcommand{\TP}{{\rm TrPic}}
\newcommand{\Pic}{{\rm Pic}}

\newcommand{\add}{{\rm add}}
\newcommand{\Tot}{{\rm Tot}}

\newcommand{\thick}{{\rm thick}}

\newcommand{\LT}{{\rm\bf L}T}
\newcommand{\AAA}{\mathcal{A}}
\newcommand{\BB}{\mathcal{B}}
\newcommand{\CC}{\mathcal{C}}
\newcommand{\DD}{\mathcal{D}}

\newcommand{\XX}{\mathcal{X}}
\newcommand{\YY}{\mathcal{Y}}
\newcommand{\KK}{\mathcal{K}}
\newcommand{\NN}{\mathcal{N}}

\newcommand{\PP}{\mathcal{P}}
\newcommand{\QQ}{\mathcal{Q}}
\newcommand{\II}{\mathcal{I}}

\newcommand{\Kb}{\mathcal{K}^{\rm b}_p}
\newcommand{\KKb}{\mathcal{K}^{\rm b}}
\newcommand{\CCb}{\mathcal{C}^{\rm b}}
\newcommand{\KbG}{\mathcal{K}^{\rm b}_{p,G}}
\newcommand{\la}{\langle}
\newcommand{\ra}{\rangle}
\renewcommand{\le}{\leqslant}
\renewcommand{\ge}{\geqslant}

\newcommand{\ot}{\otimes}
\newcommand{\kk}{\mathbf{k}}

\newenvironment{Proof}[1][Proof]{\begin{trivlist}
\item[\hskip \labelsep {\bfseries #1}]}{\flushright
$\Box$\end{trivlist}}

\numberwithin{equation}{section}

\begin{document}
\title{On standard derived equivalences of orbit categories}
\author{Yury Volkov\footnote{The author was partially supported by RFFI 13-01-00902.}\mbox{ } and Alexandra Zvonareva\footnote{The author was partially supported by RFFI 13-01-00902 and by the Chebyshev Laboratory  (Department of Mathematics and Mechanics, St. Petersburg State University) under RF Government grant 11.G34.31.0026.}}
\date{}
\maketitle
\begin{abstract}
Let $\kk$ be a commutative ring, $\AAA$ and $\BB$ -- two
$\kk$-linear categories with an action of a group $G$. We introduce
the notion of a standard $G$-equivalence from $\Kb\BB$ to $\Kb\AAA$.
We construct a map from the set of standard $G$-equivalences to the
set of standard equivalences from $\Kb\BB$ to $\Kb\AAA$ and a map
from the set of standard $G$-equivalences from $\Kb\BB$ to $\Kb\AAA$
to the set of standard equivalences from $\Kb(\BB/G)$ to
$\Kb(\AAA/G)$. We investigate the properties of these maps and apply
our results to the case where $\AAA=\BB=R$ is a Frobenius
$\kk$-algebra and $G$ is the cyclic group generated by its Nakayama
automorphism $\nu$. We apply this technique to obtain the generating
set of the derived Picard group of a Frobenius Nakayama algebra over
an algebraically closed field.
\end{abstract}

\section{Introduction}
Let $\AAA$ and $\BB$ be two derived equivalent categories. The
notion of a standard equivalence from $\DD\BB$ to $\DD\AAA$ was
introduced in \cite{Keller}. This notion generalizes the notion of a
standard equivalence for algebras \cite{Rickard}. We define such
standard equivalences in terms of tilting subcategories instead of
tilting complexes of bimodules. We denote by $\TP(\AAA,\BB)$ the set
of standard equivalences from $\Kb\BB$ to $\Kb\AAA$, standard
equivalences from $\DD\BB$ to $\DD\AAA$ correspond bijectively to
standard equivalences from $\Kb\BB$ to $\Kb\AAA$. In \cite{Keller}
it is proved that the composition of standard equivalences and the
inverse equivalence of a standard equivalence are again standard. In
particular, composition defines a group structure on
$\TP(\AAA)=\TP(\AAA,\AAA)$. We call this group the derived Picard
group of $\AAA$. In the case where a group $G$ acts on $\AAA$ and
$\BB$, we introduce the notion of a standard $G$-equivalence from
$\Kb\BB$ to $\Kb\AAA$. We denote by $\TP_G(\AAA,\BB)$ the set of
such equivalences. It appears that the composition of standard
$G$-equivalences is defined and it determines a group structure on
$\TP_G(\AAA)=\TP_G(\AAA,\AAA)$. We construct the maps
$\Phi_{\AAA,\BB}:\TP_G(\AAA, \BB)\rightarrow \TP(\AAA, \BB)$ and
$\Psi_{\AAA,\BB}:\TP_G(\AAA, \BB)\rightarrow \TP(\AAA/G, \BB/G)$
which  respect the composition. Here $\AAA/G$ is the orbit category
defined in \cite{Asashiba}. We investigate the properties of these
maps. We prove that $\Phi_{\AAA,\BB}$ sends a standard
$G$-equivalence to a Morita equivalence iff $\Psi_{\AAA,\BB}$ sends
this standard $G$-equivalence to a Morita equivalence. We prove a
theorem which gives a necessary and sufficient condition for an
element of $\TP(\AAA, \BB)$ to lie in the image of
$\Phi_{\AAA,\BB}$. In the case of a finite group $G$ we provide
necessary and sufficient condition for an element of $\TP(\AAA/G,
\BB/G)$ to lie in the image of $\Psi_{\AAA,\BB}$.

It was proved in \cite{Rickard} that the Nakayama functor commutes
with any standard derived equivalence. Suppose that $R$ is a
Frobenius algebra with a Nakayama automorphism $\nu$. Suppose that
$\ord\nu=n<\infty$. Then the cyclic group $G=\la\nu\ra\cong C_n$
acts on $R$ and we can define an algebra $R/G$. We prove that
the homomorphism $\Phi_R=\Phi_{R,R}:\TP_G(R)\rightarrow \TP(R)$ is
an epimorphism if for any $a\in Z(R)^*$ there is an element $b\in
Z(R)^*$ such that $a=b^n$. Moreover, we prove that $\Cok(\Phi_R)$ is
generated by the classes of elements from the Picard group of $R$ if
for any $a\in Z(R)^*$ there is an element $b\in R^*$ and an
automorphism $\sigma$ of $R$ such that $a=b\nu(b)\dots\nu^{n-1}(b)$
and $\sigma\nu\sigma^{-1}(x)=\nu(bxb^{-1})$ for all $x\in R$. We
apply these facts to find a generating set of the derived Picard
group of a Frobenius Nakayama algebra using the generating set of
the derived Picard group of a symmetric Nakayama algebra from
\cite{Zvonareva}.

\section{Standard equivalences as tensor products}
Throughout this paper $\kk$ is a fixed commutative ring. We assume
everywhere that $\AAA$ and $\BB$ are small, $\kk$-linear and
$\kk$-flat categories ($\AAA$ is called $\kk$-flat if $\AAA(x,y)$ is
a flat $\kk$-module for any $x,y\in\AAA$). We simply write $\ot$ for
$\ot_{\kk}$. In this section we recall some basic definitions and
results on standard equivalences of derived categories.

\begin{Def}{\rm Contravariant functors from $\AAA$ to $\Mod\kk$ are called {\it
$\AAA$-modules}. We denote by $\Mod\AAA$ the category of
$\AAA$-modules. An $\AAA$-module is called {\it projective} if it is
a direct summand of a direct sum of representable functors (a {\it
representable functor} is a functor isomorphic to $\AAA(-,x)$ for
some $x\in\AAA$). An $\AAA$-module is called {\it finitely
generated} if it is an epimorphic image of a finite direct sum of
representable functors.}
\end{Def}

A morphism of $\AAA$-modules $d:V\rightarrow V$ is called
differential if $d^2=0$. A $\mathbb{Z}$-graded $\AAA$-module is a
module $V$ with a decomposition $V=\oplus_{n\in\mathbb{Z}}V_n$. We
say that a morphism $f:V\rightarrow V'$ is of degree $m$ if
$f=\sum_{n\in\mathbb{Z}}f_n$ for some morphisms $f_n:V_n\rightarrow
V'_{n+m}$. We denote by $V[m]$ the module $V$ with the following
grading: $V[m]_n=V_{n+m}$.

\begin{Def}{\rm An {\it $\AAA$-complex} is a $\mathbb{Z}$-graded module $V$
with a differential $d_V:V\rightarrow V$ of degree $1$. A morphism
of $\AAA$-complexes is a morphism of $\AAA$-modules of degree 0
which commutes with differentials. We denote by $\CC\AAA$ the
category of $\AAA$-complexes. A morphism $f:V\rightarrow V'$ is
called {\it null homotopic} if $f=hd_V+d_{V'}h$ for some morphism of
modules $h:V\rightarrow V'$ of degree $-1$. We denote by $H(V,V')$
the space of all null homotopic morphisms from $V$ to $V'$. The {\it
homotopy category} of $\Mod\AAA$ is the category $\KK\AAA$ with the
same objects as $\CC\AAA$ and morphism spaces
$\KK\AAA(V,V')=\CC\AAA(V,V')/H(V,V')$. The {\it derived category} of
$\Mod\AAA$ denoted by $\DD\AAA$ is the localization of $\KK\AAA$ at
the set of quasi-isomorphisms (a morphism from $V$ to $V'$ is
called a quasi-isomorphism if it induces an isomorphism from $\Ker
d_V/\Im d_V$ to $\Ker d_{V'}/\Im d_{V'}$).}
\end{Def}

We denote by $\KK_p\AAA$ the full subcategory of $\KK\AAA$
consisting of all projective complexes. The canonical functor from
$\KK\AAA$ to $\DD\AAA$ induces an equivalence between categories
$\KK_p\AAA$ and $\DD\AAA$. We denote by $\Kb\AAA$ the full
subcategory of $\KK\AAA$ consisting of all finitely generated
projective complexes $V$ such that $V_n=0$ for large enough and
small enough $n\in\mathbb{Z}$.

\begin{Def} {\rm A full subcategory $\XX$ of $\Kb\AAA$ is called a {\it tilting
subcategory} for $\AAA$ if
\begin{itemize}
\item $\Kb\AAA(U, V[i])=0$ for all $U,V\in\XX$, $i\not=0$,
\item $\AAA(-,x)$ ($x\in\AAA$) lies in the smallest full triangulated subcategory of
$\Kb\AAA$ containing $\XX$ and closed under isomorphisms and direct
summands (we denote this subcategory of $\Kb\AAA$ by $\thick\XX$).
\end{itemize}}
\end{Def}

\begin{Def} {\rm The {\it tensor product} $\AAA\ot\BB$ of $\AAA$ and $\BB$ is a
$\kk$-linear category defined in the following way. Its objects are
pairs $(x,y)$ where $x\in\AAA$, $y\in\BB$. Its morphism spaces are
$$
(\AAA\ot\BB)((x_1,y_1),(x_2,y_2))=\AAA(x_1,x_2)\ot\BB(y_1,y_2).
$$
The composition in $\AAA\ot\BB$ is given by the formula
$$
(f_2\ot g_2)(f_1\ot g_1)=f_2f_1\ot g_2g_1,
$$
where $f_1\in\AAA(x_1,x_2)$, $f_2\in\AAA(x_2,x_3)$,
$g_1\in\BB(y_1,y_2)$, $g_2\in\BB(y_2,y_3)$ and $x_1,x_2,x_3\in\AAA$,
$y_1,y_2,y_3\in\BB$.}
\end{Def}

Let $X$ be an $\AAA\ot\BB^{\rm op}$-complex. It defines a functor
$T_X:\CC\BB\rightarrow\CC\AAA$ in the following way. If
$M\in\CC\BB$, then $(T_XM)(x)$ ($x\in\AAA$) is the cokernel of the
map
$$
\rho_{X,M}(x):\bigoplus\limits_{y,z\in\BB}\left(M(z)\ot \BB(y,z)\ot
X(x,y)\right)\rightarrow\bigoplus\limits_{y\in\BB}\left(M(y)\ot
X(x,y)\right)
$$
defined by the equality
$$\rho_{X,M}(x)(u\ot g\ot v)=M(g)(u)\ot v-u\ot X(\Id_x\ot g)(v)$$
for $u\in M(z)$, $v\in X(x,y)$, $g\in\BB(y,z)$. If
$f\in\AAA(x_1,x_2)$ ($x_1,x_2\in\AAA$), then $(T_XM)(f)$ is the
unique map such that the diagram
 $$
\xymatrix
{
\bigoplus_{y\in\BB}\left(M(y)\ot X(x_2,y)\right)\ar@{->}[r]\ar@{->}[d]_{\oplus_{y\in\BB}\left(\Id_{M(y)}\ot X(f\ot \Id_y)\right)}&(T_XM)(x_2)\ar@{->}[d]_{(T_XM)(f)}\\
\bigoplus_{y\in\BB}\left(M(y)\ot X(x_1,y)\right)\ar@{->}[r]&(T_XM)(x_1)
}
$$
commutes. Finally, for $\alpha\in\CC\BB(M,N)$ we obtain $T_X\alpha$
from the commutativity of diagrams
 $$
\xymatrix
{
\bigoplus_{y\in\BB}\left(M(y)\ot X(x,y)\right)\ar@{->}[r]\ar@{->}[d]_{\oplus_{y\in\BB}\left(\alpha_y\ot \Id_{X(x, y)}\right)}&(T_XM)(x)\ar@{->}[d]_{(T_X\alpha)_x}\\
\bigoplus_{y\in\BB}\left(N(y)\ot X(x,y)\right)\ar@{->}[r]&(T_XN)(x)
}
$$
We denote by $\LT_X$ the left derived functor of $T_X$, i.e. the composition
$$
\DD\BB\stackrel{\bf p}\rightarrow
\KK_p\BB\stackrel{T_X}\longrightarrow\mathcal{K}\AAA\stackrel{\pi}\rightarrow\DD\AAA,
$$
where ${\bf p}$ sends an object of $\DD\BB$ to its projective
resolution (see \cite[Theorem 3.1]{Keller} for details) and $\pi$ is
the canonical functor from $\mathcal{K}\AAA$ to $\DD\AAA$. If $X$ is
an $\AAA\ot\BB^{\rm op}$-complex, then for any $y\in\BB$ we can
define an $\AAA$-complex $X^y$ as follows: $X^y(x)=X(x,y)$ and
$X^y(f)=X(f\ot\Id_y)$ for objects and morphisms of $\AAA$
respectively. Then $\LT_X$ is an equivalence iff the following
conditions hold:
\begin{itemize}
\item $X^y$ is isomorphic to some object of $\Kb\AAA$ in $\KK\AAA$ for any $y\in\BB$,
\item the full subcategory of $\Kb\AAA$ consisting of objects isomorphic
to some $X^y$ ($y\in\BB$) in $\KK\AAA$ is a tilting subcategory for
$\AAA$,
\item the map $\BB(y,z)\rightarrow\mathcal{K}\AAA(X^y,X^z)$ is an
isomorphism for all $y,z\in\BB$.
\end{itemize}
Moreover, $\LT_X\cong\LT_Y$ iff $X\cong Y$ in $\DD(\AAA\ot\BB^{\rm
op})$. See \cite[6.1]{Keller} for details. If $\LT_X$ is an
equivalence, then $X$ is called a {\it tilting $\AAA\ot\BB^{\rm
op}$-complex}.

If $F:\AAA\rightarrow\BB$ is a functor, then we denote by $F(\AAA)$
the full subcategory of $\BB$ formed by objects isomorphic to some
$F(U)$ ($U\in\AAA$).

\begin{Def}{\rm We call $\theta:\BB\rightarrow \Kb\AAA$ a {\it tilting functor} if $\theta(\BB)$ is a tilting subcategory for $\AAA$ and $\theta$ induces an equivalence from $\BB$ to $\theta(\BB)$.}
\end{Def}

Categories $\AAA$ and $\BB$ are said to be {\it derived equivalent}
if the derived categories $\mathcal{D}\AAA$ and $\mathcal{D}\BB$ are
equivalent as triangulated categories. The following theorem is well
known (see \cite[9.2, Corollary]{Keller} and \cite[Theorem
4.6]{Asashiba}).

\begin{theorem} The following conditions are equivalent\\
$(1)$ $\AAA$ and $\BB$ are derived equivalent;\\
$(2)$ there is a tilting $\AAA\ot\BB^{\rm op}$-complex $X$;\\
$(3)$ there is a tilting functor $\theta:\BB\rightarrow\Kb\AAA$.
\end{theorem}

\begin{rema}\label{rem1} Note that if $X$ is a tilting $\AAA\ot\BB^{\rm op}$-complex, then
the objects of $\Kb\AAA$ isomorphic to $X^y$ ($y\in \BB$) form a
tilting subcategory which is equivalent to $\BB$. We denote the
corresponding tilting functor (which is defined modulo natural
isomorphism) by $\theta_X$. Conversely, if we have a tilting functor
$\theta:\BB\rightarrow\Kb\AAA$, then we can construct a tilting
$\AAA\ot\BB^{\rm op}$-complex $X$ such that $\theta_X\cong\theta$
(see \cite[Section 9]{Keller}).
\end{rema}

\begin{Def} {\rm We call an equivalence $F:\DD\BB\rightarrow\DD\AAA$ {\it standard}
if there is some $\AAA\ot\BB^{\rm op}$-complex $X$ such that
$F\cong\LT_X$. We denote by $\TP(\AAA,\BB)$ the set of all standard
equivalences from $\DD\BB$ to $\DD\AAA$ modulo natural isomorphisms.}
\end{Def}

Let $\LT_X$ be a standard equivalence. Define a $\BB\ot\AAA   ^{\rm
op}$-complex $X^T$ as follows:
$$
X^T(y,x)=\Mod\AAA(X^y,\AAA(-,x)),\,X^T(g\ot
f)=\Mod\AAA\left(X(\Id_{-}\ot g), \AAA(-,f)\right)
$$
for $x\in\AAA$, $y\in\BB$, a morphism $f$ in $\AAA$ and a morphism
$g$ in $\BB$. By \cite[6.2, Lemma]{Keller} the functor $\LT_{X^T}$
is quasi-inverse to $\LT_X$. Moreover, if $\LT_X$ and $\LT_Y$ are
standard equivalences (which can be composed), then by \cite[6.3,
Lemma]{Keller} we have $\LT_X\LT_Y\cong\LT_Z$, where $Z=T_{{\bf
p}X}Y$ and ${\bf p}X$ is the projective resolution of $X$ over
$\AAA\ot\BB^{\rm op}$.

\begin{Def}{\rm The {\it derived Picard group} of $\AAA$ is the set
$$
\TP(\AAA):=\TP(\AAA,\AAA)
$$
with the operation of composition. It follows from the arguments
above that it is actually a group.}
\end{Def}

Note that the functor $\LT_X:\DD\BB\rightarrow\DD\AAA$ is an
equivalence iff it's restriction to $\Kb\BB$ induces an equivalence
to $\Kb\AAA$. Moreover, $\LT_X\cong\LT_Y$ iff the corresponding
equivalences from $\Kb\BB$ to $\Kb\AAA$ are isomorphic. So we denote
by $\LT_X$ the corresponding equivalence from $\Kb\BB$ to $\Kb\AAA$
too. From here on we consider only standard equivalences and
identify the set $\TP(\AAA,\BB)$ with the set of standard
equivalences from $\Kb\BB$ to $\Kb\AAA$ modulo natural isomorphisms.

\section{Standard equivalences and tilting subcategories}\label{se_as_ts}

For a subcategory $E$ of $\Kb\AAA$ we denote by $\add E$ the full
subcategory of $\Kb\AAA$ consisting of direct summands of finite
direct sums of copies of objects of $E$. Let us define a category
$\CCb\add E$. Objects of $\CCb\add E$ are objects of $\add E$ with a
decomposition $V=\oplus_{n\in\mathbb{Z}}V_n$ and a differential
$d_V=\sum_{n\in \mathbb{Z}}d_{V,n}$ where $V_n\in\add E$,
$d_{V,n}\in\KK\AAA(V_n, V_{n+1})$ and $V_n=0$ for large enough and
small enough $n$. If $V,V'\in\CCb\add E$, then the set $\CCb\add
E(V,V')$ is formed by maps $f=\sum_{n\in \mathbb{Z}}f_n$ such that
$f_n\in\KK\AAA(V_n,V_n')$ and $fd_V-d_{V'}f$ equals 0 in $\KK\AAA$.
A morphism $f\in\CCb\add E(V,V')$ is called null homotopic if
$f=hd_V+d_{V'}h$ for some $h=\sum_{n\in \mathbb{Z}}h_n$ where
$h_n\in\KK\AAA(V_n,V'_{n-1})$. We denote the set of null homotopic
morphisms from $V$ to $V'$ by $B(V,V')$ again. Then $\KKb\add E$ is
a category whose objects are the same as the objects of $\CCb\add E$
and whose morphism spaces are $\KKb\add E(V,V')=\CCb\add
E(V,V')/B(V,V')$.

We denote by $\YY_{\AAA}:\AAA\rightarrow \Kb\AAA$
the Yoneda embedding, i.e. $\YY_{\AAA}(x)=\AAA(-,x)$ and
$\YY_{\AAA}(f)=\AAA(-,f)$ for an object $x$ and a morphism $f$ in
$\AAA$. Let $\theta:\BB\rightarrow \Kb\AAA$ be a tilting functor. Our aim is
to define an equivalence $F_{\theta}:\Kb\BB\rightarrow \Kb\AAA$ in
such a way that $F_{\theta}\YY_{\BB}=\theta$. Denote by $\XX$ the
category $\theta(\BB)$. Let $\PP_{\BB}$ be the category of finitely
generated projective $\BB$-modules. Let us define an equivalence
$S:\PP_{\BB}\rightarrow \add\XX$. Define $S$ on direct sums of
representable functors as follows:
$$
S\big(\oplus_{i=1}^m\BB(-,x_i)\big)=\oplus_{i=1}^m\theta(x_i).
$$
If $$f=\big(\BB(-,f_{i,j})\big)_{1\le i\le m,1\le j\le
l}:\oplus_{i=1}^m\BB(-,x_i)\rightarrow \oplus_{j=1}^l\BB(-,y_j)$$ is
a morphism in $\PP_{\BB}$, then
$$S(f)=\big(\theta(f_{i,j})\big)_{1\le i\le m,1\le j\le l}:\oplus_{i=1}^m\theta(x_i)\rightarrow \oplus_{j=1}^l\theta(y_j).$$
Let us consider an arbitrary object $U\in\PP_{\BB}$. There is some
direct sum of representable functors $W_U$ such that $U$ is a direct
summand of $W_U$. Let $\iota_U:U\rightarrow W_U$ and
$\pi_U:W_U\rightarrow U$ be the corresponding direct inclusion and
projection (for convenience we assume that $W_U=U$ and
$\iota_U=\pi_U=\Id_U$ if $U$ is a direct sum of representable
functors). It follows from \cite{LeChen} that idempotents split in
$\Kb\AAA$. In particular, they split in $\add\XX$. Since
$S(\iota_U\pi_U):S(W_U)\rightarrow S(W_U)$ is an idempotent in
$\add\XX$, there is some object $X_U\in\add\XX$ and morphisms
$\iota_U':X_U\rightarrow S(W_U)$ and $\pi_U':S(W_U)\rightarrow X_U$
such that $\pi_U'\iota_U'=\Id_{X_U}$ and
$\iota_U'\pi_U'=S(\iota_U\pi_U)$. We define $S(U)=X_U$. If
$f:U\rightarrow V$ is a morphism in $\PP_{\BB}$, then we define
$S(f)$ by the formula
$$
S(f)=\pi_V'S(\iota_Vf\pi_U)\iota_U'.
$$
It is clear that $S$ is an equivalence. Then $S$ induces an
equivalence $\bar S:\Kb\BB\rightarrow \KKb\add\XX$. Note also that $S(\iota_U)=\iota_U'$ and $S(\pi_U)=\pi_U'$.

Let us now translate some results of \cite{Rickard2} from the case
of algebras to the case of categories. Since the arguments for the
case of categories are analogous to the case of algebras, we omit
most of the proofs and give only references to the corresponding
results of \cite{Rickard2}.

Let $V$ be an object of $\CCb\add\XX$. Then we can consider $V$ as a
bigraded module $V=\oplus_{i,j\in\mathbb{Z}}V_{i,j}$ (where
$V_{i,j}=(V_j)_i$) with a differential
$\tau_0=\sum_{j\in\mathbb{Z}}d_{V_j}:V\rightarrow V$ of degree
$(1,0)$ and a morphism $\tau_1:V\rightarrow V$ of degree $(0,1)$
such that $\tau_1|_{V_{i,j}}=(-1)^{i+j}d_V|_{V_{i,j}}$. Such objects
satisfy the conditions
\begin{itemize}
\item $V_{i,j}=0$ if $i$ or $j$ is large enough or small enough;
\item $\tau_0^2=0$;
\item $\tau_0\tau_1+\tau_1\tau_0=0$;
\item $\tau_1^2=0$ in $\KK\AAA$ if we consider $V$ as an object of $\add\XX$.
\end{itemize}
We write $(V,\tau_0,\tau_1)$ for such object. Note that for two
such objects $\KK\AAA(V,V'[i]))=0$, $i \neq 0$ if we consider them as objects of
$\add\XX$. A morphism from $(V,\tau_0,\tau_1)$ to
$(V',\tau_0',\tau_1')$ in $\CCb\add\XX$ is a morphism
$\alpha:V\rightarrow V'$ of degree $(0,0)$ such that
$\alpha\tau_0=\tau_0'\alpha$ and $\alpha\tau_1-\tau_1'\alpha$ is
null homotopic if we consider $V$ and $V'$ as objects of $\add\XX$.
Moreover, a morphism $\alpha:(V,\tau_0,\tau_1)\rightarrow
(V',\tau_0',\tau_1')$ is equal to 0 in $\CCb\add\XX$ if it is null
homotopic as a morphism in $\add\XX$. If
$(V,\tau_0,\tau_1)\in\CCb\add\XX$, then we can define morphisms
$\tau_i:V\rightarrow V$ of degree $(1-i,i)$ in such a way that
$$
\sum\limits_{i=0}^l\tau_i\tau_{l-i}=0
$$
for any $l\ge 0$ (see \cite[Proposition 2.6]{Rickard2}). If
$\alpha:(V,\tau_0,\tau_1)\rightarrow (V',\tau_0',\tau_1')$, then
there is a sequence of maps $\alpha_i:V\rightarrow V'$ of degree
$(-i,i)$ such that $\alpha_0=\alpha$ and
$$
\sum\limits_{i=0}^l\alpha_i\tau_{l-i}=\sum\limits_{i=0}^l\tau_i'\alpha_{l-i}
$$
for any $l\ge 0$ (see \cite[Proposition 2.7]{Rickard2}). If
$(V,\tau_0,\tau_1)$ is an object of $\CCb\add\XX$, then we define
$\Tot (V,\tau_0,\tau_1)\in\CC\AAA$ by the formulas
$$
\Tot (V,\tau_0,\tau_1)_n=\bigoplus\limits_{i+j=n}V_{i,j},\,\,d_{\Tot
(V,\tau_0,\tau_1)}=\sum\limits_{i\ge 0}\tau_i.
$$
If $\alpha:(V,\tau_0,\tau_1)\rightarrow (V',\tau_0',\tau_1')$, then
we define $\Tot\alpha:\Tot (V,\tau_0,\tau_1)\rightarrow\Tot
(V',\tau_0',\tau_1')$ by the formula
$$
\Tot\alpha=\sum\limits_{i\ge 0}\alpha_i.
$$
Thus we define a functor $\Tot:\CCb\add\XX\rightarrow \Kb\AAA$ (see
\cite[Proposition 2.10]{Rickard2}). By \cite[Proposition
2.11]{Rickard2} the functor $\Tot$ factors through some functor $\bar
Q:\KKb\add\XX\rightarrow\Kb\AAA$. We define $F_{\theta}$ as the
composition
\begin{equation}\label{se}
\Kb\BB\stackrel{\bar S}\rightarrow\KKb\add\XX\stackrel{\bar
Q}\rightarrow\Kb\AAA.
\end{equation}
Note that $F_{\theta}$ is defined modulo isomorphism. We fix some representative of this equivalence for each tilting functor.

\begin{prop}\label{tiltsub_stand}
Let $\theta:\BB\rightarrow \Kb\AAA$ be a tilting functor and $X$ be
a tilting $\AAA\ot\BB^{\rm op}$-complex. If $\theta_X\cong\theta$,
then $\LT_X\cong F_{\theta}$.
\end{prop}
\begin{Proof} Let $F_{\theta}$ be defined by the composition \eqref{se}.
Then it is enough to prove that $T_X|_{\Kb\BB}\cong \iota\bar
Q\bar S$, where $\iota$ is the canonical embedding of $\Kb\AAA$
to $\KK\AAA$. We know that there is an isomorphism $\xi:T_X\YY_{\BB}\cong\iota\theta$. Let define an isomorphism $\zeta:T_X|_{\Kb\BB}\cong \iota\bar
Q\bar S$. Let $U=\oplus_{i=1}^m\BB(-,x_i)$ be a direct sum of representable functors. Then we define
$$\zeta_U:T_X(U)=\oplus_{i=1}^mT_X\YY_{\BB}(x_i)\rightarrow \oplus_{i=1}^m\iota\theta(x_i)=\iota\bar Q\bar S(U)$$
by the equality $\zeta_U=\oplus_{i=1}^m\xi_{x_i}$. For $U\in\PP_{\BB}$ we define
$$
\zeta_U=S(\pi_U)\zeta_{W_U}T_X(\iota_U):T_X(U)\rightarrow S(U)=\iota\bar Q\bar S(U).
$$
It is easy to see that $\zeta$ defines an isomorphism from $T_X|_{\PP_{\BB}}$ to $\iota\bar Q\bar S|_{\PP_{\BB}}$.

Let now $U=\oplus_{n\in\mathbb{Z}}U_n$ be an arbitrary object of $\Kb\BB$. Then $T_X(U)$ is a totalization of a bigraded module $V=\oplus_{i,j\in\mathbb{Z}}T_X(U_j)_i$ with differential $\tau_0=\sum_{j\in\mathbb{Z}}d_{T_X(U_j)}$ of degree $(1,0)$ and differential $\tau_1$ of degree $(0,1)$ defined by the equality $\tau_1|_{V_{i,j}}=(-1)^{i+j}T_X(d_{U})|_{T_X(U_j)_i}$. At the same time $\iota\bar Q\bar S(U)$ is a totalization of a bigraded module $V'=\oplus_{i,j\in\mathbb{Z}}S(U_j)_i$ with some differentials $\tau'_i$ ($i\ge 0$) of degree $(1-i, i)$ such that
$\tau'_0=\sum_{j\in\mathbb{Z}}d_{S(U_j)}$ and $\tau'_1|_{V_{i,j}}=(-1)^{i+j}S(d_U)|_{S(U_j)_i}$. Here we write $S(d_U)$ for some representative of homotopy class of it. If we choose a representative of homotopy class of $\zeta_{U_j}$ for all $j\in\mathbb{Z}$, then we obtain a differential $\zeta_0:V\rightarrow V'$ of degree $(0,0)$ such that $\tau'_0\zeta_0=\zeta_0\tau_0$ and $\tau'_1\zeta_0$ is homotopic to $\zeta_0\tau_1$ if we consider $V$ and $V'$ as objects of $\Kb\AAA$ (i.e. if we forget the grading on $U$). Analogously to \cite[Proposition 2.7]{Rickard2} we can construct $\zeta_i$ for $i\ge 0$ of degree $(-i,i)$ such that
$$
\zeta_l\tau_0+\zeta_{l-1}\tau_1=\sum\limits_{i=0}^l\tau_i'\zeta_{l-i}
$$
for any $l\ge 1$. We define $\zeta_U=\sum\limits_{i\ge 0}\zeta_i:\Tot(V,\tau_0,\tau_1)\rightarrow\Tot(V',\tau_0',\tau_1')$. If $U,U'\in \Kb\BB$ and $f\in\CCb\BB(U,U')$, then it is clear that $\zeta_{U'}T_X(f)-\bar Q\bar S(f)\zeta_U$ equals $\sum\limits_{i\ge 0}\upsilon_i$, where $\upsilon_i$ is of degree $(-i,i)$ and $\upsilon_0$ is null homotopic. Then it follows from arguments above that $\zeta:T_X|_{\Kb\BB}\rightarrow \iota\bar Q\bar S$ is a morphism of functors. It is clear that it is actually an isomorphism.

\end{Proof}

Note that by Proposition \ref{tiltsub_stand} and Remark \ref{rem1}
an equivalence $F:\Kb\BB\rightarrow\Kb\AAA$ is standard iff $F\cong
F_{\theta}$ for some tilting functor $\theta:\BB\rightarrow\Kb\AAA$.
Moreover, $F_{\theta}\cong F_{\theta'}$ iff $\theta\cong \theta'$.

\section{$G$-functors and orbit categories}
We say that $G$ acts on the category $\AAA$ if there is a
homomorphism of groups $\Delta:G\rightarrow\Aut(\AAA)$. In this case
we simply write $g$ instead of $\Delta(g)$ for $g\in G$. Throughout
this section we assume that $\AAA$ and $\BB$ are categories with
$G$-action. Now we recall some definitions from \cite{Asashiba}.

\begin{Def} {\rm A family $\eta=(\eta_g)_{g\in G}$ of natural isomorphisms
$\eta_g:F\circ g\rightarrow g\circ F$ is called a {\it
$G$-equivariance adjuster} for the functor $F:\AAA\rightarrow\BB$ if
the diagram
$$
\xymatrix
{
(F\circ gh)x\ar@{->}[r]^{\eta_{g,hx}}\ar@{->}[rd]_{\eta_{gh,x}}&(g\circ F\circ h)x\ar@{->}[d]^{g(\eta_{h,x})}\\
&(gh\circ F)x\\
}
$$
commutes for all $g,h\in G$ and $x\in \AAA$. We say that $F$ is a
{\it $G$-equivariant functor} if there is a $G$-equivariance
adjuster for $F$. The functor $F$ is called {\it strictly
$G$-equivariant} if $F\circ g$ equals $g\circ F$, i.e. if the
$G$-equivariance adjuster can be set to be identity.}
\end{Def}

\begin{Def}{\rm
A {\it $G$-functor} from $\AAA'$ to $\AAA$ is a pair $(F,\eta)$,
where $F$ is a functor from $\AAA'$ to $\AAA$ and $\eta$ is a
$G$-equivariance adjuster for $F$. A morphism from $(F,\eta)$ to
$(F',\eta')$ is a morphism of functors $\alpha:F\rightarrow F'$ such
that the diagram
 $$
\xymatrix
{
Fg\ar@{->}[r]^{\alpha g}\ar@{->}[d]_{\eta_g}&F'g\ar@{->}[d]_{\eta'_g}\\
gF\ar@{->}[r]^{g\alpha}&gF'
}
$$
commutes for any $g\in G$. It is clear that a morphism of
$G$-functors is an isomorphism iff it is an isomorphism of functors.
A $G$-functor $(F,\eta)$ is called a {\it $G$-equivalence} if $F$ is
an equivalence. If $F$ is strictly $G$-equivariant, we simply write
$F$ for the corresponding $G$-functor.}
\end{Def}

If $(F,\eta):\AAA'\rightarrow \AAA$ and
$(F',\eta'):\AAA''\rightarrow\AAA'$ are $G$-functors, then we define
their composition by the formula
$$
(F,\eta)\circ(F',\eta')=(FF',\eta F(\eta')):\AAA''\rightarrow \AAA,
$$
where $$\big(\eta F(\eta')\big)_{g,x}=\eta_{g,F'x}\circ
F(\eta'_{g,x}):FF'gx\rightarrow gFF'x.$$ It is easy to see that the
composition defined above is associative (see \cite[Lemma
2.8]{Asashiba2}). Moreover, it respects isomorphisms of
$G$-functors. If $(F,\eta):\AAA'\rightarrow\AAA$ is a
$G$-equivalence and $\bar F$ is an equivalence quasi-inverse to $F$,
then there is $\bar\eta=(\bar\eta_g)_{g\in G}$ ($\bar\eta_g:\bar
Fg\rightarrow g\bar F$) such that $(\bar
F,\bar\eta):\AAA\rightarrow\AAA'$ is a $G$-equivalence and
$$
(F,\eta)\circ(\bar F,\bar\eta)\cong \Id_{\AAA}\mbox{ and }(\bar F,\bar\eta)\circ(F,\eta)\cong\Id_{\AAA'}.
$$
This follows from the proof of \cite[Theorem 9.1]{Asashiba2}. We
call this $G$-equivalence $(\bar F,\bar\eta)$ (which is defined by
$(F,\eta)$ modulo isomorphism of $G$-equivalences) the quasi-inverse
$G$-equivalence to $(F,\eta)$.

\begin{rema}\label{Gfun_ind}
Let $F,F':\AAA\rightarrow\AAA'$ be functors, $\eta$ -- a
$G$-equivariance adjuster for $F$ and $\xi:F\rightarrow F'$ -- an
isomorphism. Then $\eta'=(\eta'_g)_{g\in G}$, where
$\eta'_{g,x}=g(\xi_x)\circ\eta_{g,x}\circ \xi^{-1}_{gx}$, is a
$G$-equivariance adjuster for $F'$. Moreover,
$(F',\eta')\cong(F,\eta)$.
\end{rema}

\begin{Def}{ \rm The {\it orbit category} $\AAA/G$ is defined as follows.
\begin{itemize}
\item The class of objects of $\AAA/G$ is equal to that of $\AAA$.
\item Let $x,y\in\AAA/G$. The set $\AAA/G(x, y)$ consists of
$f=(f_{h,g})_{g,h\in G}$ such that
\begin{itemize}
\item $f_{h,g}\in\AAA(gx,hy)$;
\item the sets $\{g\in G\mid f_{g,h}\not=0\}$ and $\{g\in G\mid f_{h,g}\not=0\}$ are finite for any $h\in G$;
\item $f_{lh,lg}=l(f_{h,g})$ for all $g,h,l\in G$.
\end{itemize}
\item The composition in $\AAA/G$ is defined by the fomula
$$
(f'_{h,g})_{g,h\in G}(f_{h,g})_{g,h\in G}=\left(\sum\limits_{l\in
G}f'_{h,l}f_{l,g}\right)_{g,h\in G}.
$$
\end{itemize}}
\end{Def}

We can define the action of $G$ on the category $\Mod\AAA$ by the
formula ${}^gX:=X\circ g^{-1}$ for $X\in\Mod\AAA$ and in the obvious
way for morphisms. Note that ${}^g\AAA(-,x)\cong\AAA(-,g(x))$. This
action of $G$ induces an action of $G$ on the category $\Kb\AAA$.
Let $\theta:\BB\rightarrow \Kb\AAA$ be a tilting functor. If $\theta$ is
$G$-equivariant, then the categories $\AAA/G$ and $\BB/G$ are derived equivalent by
\cite[Theorem 4.11]{Asashiba}.

\begin{rema}
Let $X$ be a tilting $\AAA\ot\BB^{\rm op}$-complex. The category
$\AAA\ot\BB^{\rm op}$ can be equipped with the diagonal action of
$G$, i.e. for $g \in G$ we put $g(x,y)=(gx,gy)$ and $g(f\ot f')=gf\ot gf'$ for
$x\in\AAA$, $y\in\BB$ and morphisms $f$ in $\AAA$, $f'$ in $\BB$.
Then it is easy to see that a $G$-equivariance adjuster for
$\theta_X:\BB\rightarrow\Kb\AAA$ is the same thing as a family of maps
$\phi_g\in\CC(\AAA\ot\BB^{\rm op})(X, {}^gX)$ such that the diagram
$$
\xymatrix
{
X^y\ar@{->}[r]^{(\phi_{g})_y}\ar@{->}[rd]_{(\phi_{gh})_y}&({}^gX)^y\ar@{->}[d]^{({}^g\phi_{h})_y}\\
&({}^{gh}X)^y\\
}
$$
commutes in $\KK\AAA$ for all $g,h\in G$ and $y\in\BB$.
\end{rema}

\section{Standard $G$-equivalences}

From here on we equip the category $\AAA/G$ with the trivial action
of $G$ (i.e. $\Delta(g)=\Id_{\AAA/G}$ for all $g\in G$) for any
category $\AAA$ with a $G$-action. We equip the category
$\Kb(\AAA/G)$ with the trivial action of $G$ as well.

\begin{Def}{\rm The {\it canonical functor} $P:\AAA\rightarrow \AAA/G$ is defined by
$P(x)=x$ and $P(f)=\big(\delta_{g,h}g(f)\big)_{g,h\in G}$ for
$x,y\in\AAA$ and $f\in\AAA(x,y)$. Let $s=(s_g)_{g\in G}$ be the
collection of maps $s_g:Pg\rightarrow P$, where
$s_{g,x}=(\delta_{hg,h'}\Id_{h'x})_{h,h'\in G}:Pgx\rightarrow Px$.
Then $s$ is a $G$-equivariance adjuster for $P$. We call $(P,s)$ the
{\it canonical $G$-functor}. }
\end{Def}

By \cite[Proposition 2.6]{Asashiba} every morphism in
$f\in\AAA/G(x,y)$ can be uniquely presented in the form
\begin{equation}\label{canon_form}
f=\sum\limits_{g\in G}s_{g,y}\circ Pf_g
\end{equation}
for some $f_g\in\AAA(x,gy)$.

\begin{Def}{\rm We define the {\it pullup functor}
$P^{\bullet}:\Mod\AAA/G\rightarrow\Mod\AAA$ by the formula
$P^{\bullet}(X)=X\circ P$ for all $X\in \Mod\AAA/G$. The {\it
pushdown functor} $P_{\bullet}:\Mod\AAA\rightarrow\Mod\AAA/G$ is the
functor left adjoint to $P^{\bullet}$. It also induces a functor
$P_{\bullet}:\Kb\AAA\rightarrow\Kb(\AAA/G)$. }
\end{Def}

We will use the explicit description of $P_{\bullet}$ obtained in
\cite[Theorem 4.3]{Asashiba}. In particular, we have
$(P_{\bullet}X)(x)=\oplus_{g\in G}X(gx)$. The same theorem says that
the map $s_{\bullet}$ defined by the commutative diagram
 $$
 \xymatrix
 {
 (P_{\bullet}{}^gX)(x)\ar@{=}[d]\ar@{->}[rrrr]^{s_{\bullet,g,X,x}}&&&&(P_{\bullet}X)(x)\ar@{=}[d]\\
 \bigoplus\limits_{h\in G}X(g^{-1}hx)\ar@{->}[rrrr]^{(\delta_{g^{-1}h,h'}\Id_{X(h'x)})_{h,h'\in G}}&&&&\bigoplus\limits_{h'\in G}X(h'x)
 }
 $$
is a $G$-equivariance adjuster for $P_{\bullet}$. Moreover, by
\cite[Theorem 4.4]{Asashiba} every morphism
$f\in\Kb(\AAA/G)(P_{\bullet}X,P_{\bullet}Y)$ can be uniquely
presented in the form
\begin{equation}\label{P_dot_rep}
f=\sum\limits_{g\in G}s_{\bullet,g,Y}\circ P_{\bullet}f_g
\end{equation}
for some $f_g\in \Kb\AAA(X,{}^gY)$. In addition, we have an
isomorphism $\gamma_y:\BB/G(-,Py)\rightarrow P_{\bullet}\BB(-,y)$
($y\in\BB$) defined by the formula $\gamma_y(s_{g,y}Pf)=g^{-1}(f)\in
\BB(g^{-1}x,y)$ for $f\in\BB(x, gy)$.

Note also that the Yoneda embedding $\YY_{\AAA}:\AAA\rightarrow
\Kb\AAA$ admits a $G$-equivariance adjuster $\phi=(\phi_g)_{g\in G}$
defined by the formula $\phi_{g,x}(f)=g^{-1}(f)$ for
$f\in\AAA(y,gx)$.

\begin{Lemma}
For all $g\in G$, $x,y\in\AAA$ and $f\in\AAA(x,gy)$ the following equality holds
\begin{equation}\label{phi_gamma}
\gamma_y^{-1}s_{\bullet,g,\AAA(-,y)}P_{\bullet}\phi_{g,y}P_{\bullet}\YY_{\AAA}(f)\gamma_x=\YY_{\AAA/G}(s_{g,y}Pf).
\end{equation}
\end{Lemma}
\begin{Proof}
It is enough to prove that the left and the right parts of the equality \eqref{phi_gamma} send the element $s_{h,x}Pf'\in\AAA/G(Pz,Px)$ to the same element of $\AAA/G(Pz,Py)$ for all $h\in G$, $z\in\AAA$ and $f'\in\AAA(z,hx)$. Direct calculations show that both parts of \eqref{phi_gamma} send $s_{h,x}Pf'$ to $s_{hg,y}P(h(f)f')$.
\end{Proof}

\begin{Def}
{\rm A $G$-equivalence $(F,\eta):\Kb \BB\rightarrow \Kb \AAA$ is called {\it standard $G$-equivalence} if $F$ is a standard equivalence and there is a standard equivalence $F':\Kb(\BB/G)\rightarrow \Kb(\AAA/G)$ such that there is an isomorphism of $G$-functors
\begin{equation}\label{standard_G}
(P_{\bullet},s_{\bullet})\circ (F,\eta)\cong F'(P_{\bullet},s_{\bullet}).
\end{equation}
We denote by $\TP_G(\AAA,\BB)$ the set of isomorphism classes of
standard $G$-equivalences from $\Kb \BB$ to $\Kb \AAA$. }
\end{Def}

It is clear that the composition of standard $G$-equivalences and
the quasi-inverse $G$-equivalence to a standard $G$-equivalence are
standard.

\begin{Def}{\rm The {\it derived Picard $G$-group} of $\AAA$ is the set
$$
\TP_G(\AAA):=\TP_G(\AAA,\AAA)
$$
with the operation of composition. It follows from the arguments
above that it is actually a group.}
\end{Def}

Let $\theta:\BB\rightarrow \Kb\AAA$ be a tilting functor and $\psi$
be a $G$-equivariance adjuster for $\theta$. We denote $\theta(\BB)$
by $\XX$, note that by agreement $\theta(\BB)$ is closed under isomorphism. Then the category $P_{\bullet}\XX$ is a tilting
subcategory for $\AAA/G$ (see the proof of \cite[Theorem
4.7]{Asashiba}). We will construct a tilting functor
$\mu_{\theta,\psi}:\BB/G\rightarrow\Kb(\AAA/G)$ which induces
an equivalence from $\BB/G$ to $P_{\bullet}\XX$. We define it on
objects by the formula $\mu_{\theta,\psi}(Py)=P_{\bullet}\theta(y)$
for $y\in\BB$. Let us consider a morphism $f=\sum\limits_{g\in
G}s_{g,y} Pf_g\in\BB/G(Px,Py)$. Define
$$\mu_{\theta,\psi}(f)=\sum\limits_{g\in G}s_{\bullet,g,\theta(y)}\circ P_{\bullet}\psi_{g,y}\circ P_{\bullet}\theta(f_g).$$
It is easy to check that $\mu_{\theta,\psi}$ is a functor. Since
$\theta$ induces an equivalence to $\XX$, $\mu_{\theta,\psi}$
induces an equivalence to $P_{\bullet}\XX$ by the arguments above.
So an equivalence
$F_{\mu_{\theta,\psi}}:\Kb(\BB/G)\rightarrow\Kb(\AAA/G)$ is defined.
Note that if $(\theta,\psi)\cong(\theta',\psi')$, then
$\mu_{\theta,\psi}\cong \mu_{\theta',\psi'}$.

\begin{prop}\label{sur_in} There is a $G$-equivariance adjuster $\eta$ for $F_{\theta}$ such
that
$$F_{\mu_{\theta,\psi}}\circ(P_{\bullet},s_{\bullet})\cong(P_{\bullet},s_{\bullet})\circ (F_{\theta},\eta).$$
\end{prop}
\begin{Proof} Note that ${}^gU$ lies in $\XX$ for all $g\in G$ and $U\in\XX$.
Indeed, since $\theta$ induces an equivalence to $\XX$, there is
some $y\in\BB$ such that $\theta(y)\cong U$. Then ${}^gU\cong
{}^g\theta(y)\cong \theta(gy)$ and ${}^gU$ lies in $\XX$ because
$\XX$ is closed under isomorphisms. Then it is easy to see that the
action of $G$ on $\Kb\AAA$ induces an action on $\KKb\add\XX$ and
the $G$-functor
$(P_{\bullet},s_{\bullet}):\Kb\AAA\rightarrow\Kb(\AAA/G)$ induces a
$G$-functor $(P_{\bullet},s_{\bullet}):\KKb\add \XX\rightarrow
\KKb\add P_{\bullet}\XX$ (here we equip the category $\KKb\add
P_{\bullet}\XX$ with the trivial action of $G$). Let us consider the
diagram
\begin{equation}\label{commute}
\xymatrix
{
\Kb\BB\ar@{->}[rr]^{(\bar S,\eta_S)}\ar@{->}[d]^{(P_{\bullet},s_{\bullet})}&&\KKb\add\XX\ar@{->}[rr]^{(\bar Q,\eta_Q)}\ar@{->}[d]^{(P_{\bullet},s_{\bullet})}&&\Kb\AAA\ar@{->}[d]^{(P_{\bullet},s_{\bullet})}\\
\Kb(\BB/G)\ar@{->}[rr]^{\bar S_G}&&\KKb\add P_{\bullet}\XX\ar@{->}[rr]^{\bar Q_G}&&\Kb(\AAA/G)
}
\end{equation}
where the rows are the compositions corresponding to \eqref{se} from
the construction of $F_{\theta}$ and $F_{\mu_{\theta,\psi}}$ (if we
omit the $G$-equivariance adjusters in the upper row). It is enough
to show that $\eta_S$ and $\eta_Q$ can be constructed in such a way
that the diagram \eqref{commute} becomes commutative modulo
isomorphism as a diagram of $G$-functors. Here we consider the
functors in the lower row as strict $G$-functors.

It is clear that $\Kb\BB=\KKb\add\YY_{\BB}(\BB)$. Let us define a
$G$-equivariance adjuster $\eta$ for the functor
$S:\add\YY_{\BB}(\BB)=\PP_{\BB}\rightarrow\add\XX$ (see section
\ref{se_as_ts}) in the following way. If
$U=\oplus_{i=1}^m\BB(-,x_i)$ is a direct sum of representable
functors, then
$$
\eta_{g,U}=\left(\oplus_{i=1}^m\psi_{g,x_i}\right)S(\oplus_{i=1}^m\phi_{g,x_i}^{-1}):S({}^gU)\cong\oplus_{i=1}^m\theta(gx_i)\rightarrow \oplus_{i=1}^m{}^g\theta(x_i)={}^gS(U).
$$
Let us now consider an arbitrary $U\in \PP_{\BB}$. Then we define
$\eta_{g,U}$ by the formula
$$
\eta_{g,U}={}^gS(\pi_U)\eta_{g,W_U}S({}^g\iota_U)
$$
(see the construction of the functor $S$ for notation). Direct calculations involving formula \eqref{phi_gamma} show that
$(S,\eta)$ is a $G$-functor. Then $\eta$ induces a
$G$-equivariance adjuster $\eta_S$ for $\bar S$ in the obvious way.
To prove the commutativity of the first square in \eqref{commute} it
is enough to prove that the diagram
$$
\xymatrix
{
\PP_{\BB}\ar@{->}[rr]^{(S,\eta)}\ar@{->}[d]^{(P_{\bullet},s_{\bullet})}&&\add\XX\ar@{->}[d]^{(P_{\bullet},s_{\bullet})}\\
\PP_{\BB/G}\ar@{->}[rr]^{S_G}&&\add P_{\bullet}\XX
}
$$
commutes modulo isomorphism of $G$-functors. Let us construct an
isomorphism $\chi:P_{\bullet}S\rightarrow S_GP_{\bullet}$. If
$U=\oplus_{i=1}^m\BB(-,x_i)$ is a direct sum of representable
functors, then
$$P_{\bullet}S(U)=\oplus_{i=1}^mP_{\bullet}\theta(x_i)=S_G(\oplus_{i=1}^m\BB/G(-,Px_i)).$$
We set $\chi_U=S_G(\oplus_{i=1}^m\gamma_{x_i})$. For an arbitrary
$U\in \PP_{\BB}$ we define $\chi_U$ by the formula
$$
\chi_U=S_G(P_{\bullet}\pi_U)\chi_{W_U}P_{\bullet}S(\iota_U).
$$
Direct calculations involving formula \eqref{phi_gamma} show that $\chi$ is the required isomorphism of
$G$-functors.

It remains to check the commutativity of the second square in
\eqref{commute}. Let us take an object of $\KKb\add\XX$ represented
by a triple $(U,\tau_0,\tau_1)$ (see section \ref{se_as_ts}). Then
$P_{\bullet}U$ can be represented by the triple
$(P_{\bullet}U,P_{\bullet}\tau_0,P_{\bullet}\tau_1)$. Suppose that
$\bar Q$ sends $(U,\tau_0,\tau_1)$ to the totalization
$(U,\sum_{i\ge 0}\tau_i)$ and  $\bar Q_G$ sends
$(P_{\bullet}U,P_{\bullet}\tau_0,P_{\bullet}\tau_1)$ to the
totalization $(P_{\bullet}U,\sum_{i\ge 0}\upsilon_i)$. It is clear
that $\upsilon_0=P_{\bullet}\tau_0$ and that $\upsilon_1$ is
homotopic to $P_{\bullet}\tau_1$ if we consider $P_{\bullet}U$ as an
object of $\Kb(\AAA/G)$. By the results of section \ref{se_as_ts}
there is a sequence of $\AAA/G$-module morphisms
$\alpha_i:P_{\bullet}U\rightarrow P_{\bullet}U$ ($i\ge 0$) such that
\begin{itemize}
\item $\alpha_i$ is of degree $(-i,i)$,
\item $\alpha_0=\Id_{P_{\bullet}U}$,
\item $\sum\limits_{i=0}^l\alpha_iP_{\bullet}\tau_{l-i}=\sum\limits_{i=0}^l\upsilon_i\alpha_{l-i}$.
\end{itemize}
In this case define the isomorphism $\zeta_U$ from the totalization
$(P_{\bullet}U,\sum_{i\ge 0}P_{\bullet}\tau_i)$ to the totalization
$(P_{\bullet}U,\sum_{i\ge 0}\upsilon_i)$ by the formula
$\zeta_U=\sum_{i\ge 0}\alpha_i$. Let $U, V\in \KKb\add\XX$,
$f:U\rightarrow V$ be a morphism  in $\CCb\add\XX$. It is clear that
the map $\zeta_VP_{\bullet}\bar Q(f)-\bar Q_GP_{\bullet}(f)\zeta_U$
is a totalization of a map from $P_{\bullet}U$ to $P_{\bullet}V$
which have nonzero components only in degrees $(-i,i)$ for $i>0$. It
follows from the results of section \ref{se_as_ts} that the
totalization of such a map is null homotopic. So $\zeta$ gives an
isomorphism from $P_{\bullet}\bar Q$ to $\bar Q_GP_{\bullet}$. It
remains to construct a $G$-equivariance adjuster $\eta_Q$ for $\bar
Q$ such that the diagram
 $$
\xymatrix
{
P_{\bullet}\bar Q{}^gU\ar@{->}[r]^{\zeta_{{}^gU}}\ar@{->}[d]_{s_{\bullet,g,\bar QU}\circ P_{\bullet}\eta_{Q,g,U}}&\bar Q_GP_{\bullet}{}^gU\ar@{->}[d]_{\bar Q_G(s_{\bullet,g,U})}\\
P_{\bullet}\bar QU\ar@{->}[r]^{\zeta_U}& \bar Q_GP_{\bullet}U
}
$$
commutes for all $U\in \KKb\add\XX$ and $g\in G$. The construction
of such isomorphisms $\eta_{Q,g,U}:\bar Q{}^gU\rightarrow {}^g\bar
QU$ is analogous to the construction of $\zeta_U$ and so it is left
to the reader.
\end{Proof}

\begin{coro}
Let $\theta:\BB\rightarrow \Kb\AAA$ be a tilting functor. Then the following statements are equivalent:\\
$1)$ there is a $G$-equivariance adjuster for $\theta$;\\
$2)$ there is a $G$-equivariance adjuster for $F_{\theta}$;\\
$3)$ there is a $G$-equivariance adjuster $\eta$ for $F_{\theta}$ such that $(F_{\theta},\eta)$ is a standard $G$-equivalence.
\end{coro}
\begin{Proof}
The implication $"1)\Rightarrow 3)"$ follows from Proposition \ref{sur_in}. Implications $"3)\Rightarrow 2)\Rightarrow 1)"$ are obvious.
\end{Proof}

Let $(F_{\theta},\eta):\Kb\BB\rightarrow\Kb\AAA$ be a $G$-equivalence. We define $\psi_{g,y}:\theta(gy)\rightarrow{}^g\theta(y)$ by the formula
$$
\psi_{g,y}=\eta_{g,B(-,y)}F_{\theta}(\phi_{g,y}).
$$
It is clear that $\psi=(\psi_g)_{g\in G}$ is a $G$-equivariance adjuster for $\theta$.

\begin{theorem}\label{Psi_cor}
Let $(F_{\theta},\eta):\Kb\BB\rightarrow\Kb\AAA$ be a standard
$G$-equivalence. Then $F_{\mu_{\theta,\psi}}$ is determined by the
condition \eqref{standard_G} uniquely modulo isomorphism.
\end{theorem}
\begin{Proof} Since $(F_{\theta},\eta)$ is a standard $G$-equivalence, there is some standard equivalence $F':\Kb(\BB/G)\rightarrow\Kb(\AAA/G)$ satisfying the condition \eqref{standard_G}.

It is enough to prove that $F'\YY_{\BB/G}\cong\mu_{\theta,\psi}$.
From \eqref{standard_G} we have a natural isomorphism
$\xi:F'P_{\bullet}\rightarrow P_{\bullet}F_{\theta}$ such that the
diagram
\begin{equation}\label{G_fun}
\begin{array}{c}
\xymatrix
{
F'P_{\bullet}{}^gU\ar@{->}[rr]^{\xi_{{}^gU}}\ar@{->}[d]_{F'(s_{\bullet,g,U})}&&P_{\bullet}F_{\theta}{}^gU\ar@{->}[d]_{s_{\bullet,g,F_{\theta}U}\circ P_{\bullet}(\eta_{g,U})}\\
F'P_{\bullet}U\ar@{->}[rr]^{\xi_U}&&P_{\bullet}F_{\theta}U
}
\end{array}
\end{equation}
commutes for any $U\in\Kb\BB$. Let us define
$\zeta_{Py}:F'\YY_{\BB/G}(Py)\rightarrow\mu_{\theta,\psi}(Py)$ by
the formula
$$\zeta_{Py}=\xi_{\BB(-,y)}F'\gamma_y.$$
Then using \eqref{phi_gamma}, \eqref{G_fun} and the fact that $\xi$ is a morphism of functors we get
$$
\begin{aligned}
&\zeta_{Py}F'\YY_{\BB/G}(s_{g,y}Pf)=\zeta_{Py}F'\big(\gamma_y^{-1}s_{\bullet,g,B(-,y)}P_{\bullet}(\phi_{g,y})\gamma_{gy}\BB/G(-,Pf)\big)\\
=&s_{\bullet,g,\theta(y)}\circ P_{\bullet}\eta_{g,\BB(-,y)}\circ \xi_{{}^g\BB(-,y)}F'P_{\bullet}(\phi_{g,y}\BB(-,f))F'(\gamma_x)=\mu_{\theta,\psi}(s_{g,y}Pf)\zeta_{Px}.
\end{aligned}
$$
for all $x,y\in\BB$, $g\in G$, $f\in\BB(x,gy)$. Since any morphism
in $\BB/G$ is of the form \eqref{canon_form}, $\zeta$ is the
required isomorphism from $F'\YY_{\BB/G}$ to $\mu_{\theta,\psi}$.
\end{Proof}

\section{The maps $\Phi$ and $\Psi$}

In this section we define two maps:
$$
\Psi_{\AAA,\BB}:\TP_G(\AAA,\BB)\rightarrow \TP(\AAA/G,\BB/G)\mbox{ and }\Phi_{\AAA,\BB}:\TP_G(\AAA,\BB)\rightarrow \TP(\AAA,\BB).
$$
Then we investigate some of their properties.

Let $(F,\eta):\Kb\BB\rightarrow \Kb\AAA$ be a standard
$G$-equivalence. Then we define $\Phi_{\AAA,\BB}$ as follows:
$$
\Phi_{\AAA,\BB}(F,\eta)=F
$$
and define $\Psi_{\AAA,\BB}(F,\eta)$ to be the unique standard
equivalence $F'$ satisfying the condition \eqref{standard_G}. The
correctness of the definition of $\Psi_{\AAA,\BB}$ follows from
Theorem \ref{Psi_cor}. It is clear that
$$
\Phi_{\AAA,\BB'}\Big((F,\eta)\circ(F',\eta')\Big)=\Phi_{\AAA,\BB}(F,\eta)\circ\Phi_{\BB,\BB'}(F',\eta')
$$
and
$$
\Psi_{\AAA,\BB'}\Big((F,\eta)\circ(F',\eta')\Big)=\Psi_{\AAA,\BB}(F,\eta)\circ\Psi_{\BB,\BB'}(F',\eta').
$$
In particular, $\Phi_{\AAA}:=\Phi_{\AAA,\AAA}$ and $\Psi_{\AAA}:=\Psi_{\AAA,\AAA}$ are homomorphisms of groups.

\begin{Def}{\rm A standard equivalence $F:\Kb\BB\rightarrow\Kb\AAA$ is called a {\it
Morita equivalence} if $F\cong F_{\theta}$ where $\theta(y)$ is
isomorphic to some object $U$ concentrated in degree 0 ($U_n=0$ for
$n\not=0$) for any $y\in\BB$. We denote by $\Pic(\AAA,\BB)$ the set
of Morita equivalences from $\BB$ to $\AAA$ modulo isomorphisms. It
is clear that the composition of Morita equivalences and the inverse
to a Morita equivalence are again Morita equivalences. In
particular, the set $\Pic(\AAA):=\Pic(\AAA,\AAA)$ is a subgroup of
$\TP(\AAA)$. This group is called the {\it Picard group} of $\AAA$.}
\end{Def}

\begin{theorem}\label{Picard} Let $(F,\eta):\Kb\BB\rightarrow\Kb\AAA$ be a standard $G$-equivalence. Then
$$
\Phi_{\AAA,\BB}(F,\eta)\in\Pic(\AAA,\BB)\Leftrightarrow\Psi_{\AAA,\BB}(F,\eta)\in\Pic(\AAA/G,\BB/G).
$$
In particular,
$$
\Phi_{\AAA}^{-1}\big(\Pic(\AAA)\big)=\Psi_{\AAA}^{-1}\big(\Pic(\AAA/G)\big).
$$
\end{theorem}
\begin{Proof} By Remark \ref{Gfun_ind} we can assume that $F=F_{\theta}$ for some
tilting functor $\theta$. Denote $\theta(\BB)$ by $\XX$. By Theorem
\ref{Psi_cor} we have $\Psi_{\AAA,\BB}(F,\eta)\cong F_{\mu}$ for
some equivalence $\mu:\BB/G\rightarrow \Kb(\AAA/G)$ such that
$\mu(Py)=P_{\bullet}\theta(y)$ for any $y\in\BB$.

Suppose that $F\in\Pic(\AAA,\BB)$. Let us consider $y\in\BB$. There
is some object $U\in\Kb\AAA$ concentrated in degree 0 such that
$\theta(y)\cong U$. Then $\mu(Py)=P_{\bullet}\theta(y)\cong
P_{\bullet}U$. It is clear that $P_{\bullet}U$ is concentrated in
degree 0. Consequently,
$\Psi_{\AAA,\BB}(F,\eta)\in\Pic(\AAA/G,\BB/G)$.

Suppose now that $\Psi_{\AAA,\BB}(F,\eta)\in\Pic(\AAA/G,\BB/G)$. It
is enough to prove that any object of $\XX$ is isomorphic in
$\Kb\AAA$ to some object concentrated in degree 0. Any object of
$P_{\bullet}\XX$ is isomorphic in $\Kb(\AAA/G)$ to some object
concentrated in degree 0 by our assumption. Consider some $U\in\XX$.
We know that $P_{\bullet}U$ is isomorphic to an object concentrated
in degree 0. Then $P^{\bullet}P_{\bullet}U$ is isomorphic to an
object concentrated in degree 0 in $\KK_p\AAA$. Since $U$ is a
direct summand of $P^{\bullet}P_{\bullet}U$ (see the proofs of
\cite[Theorems 4.3 and 4.4]{Asashiba}), $U$ is isomorphic to some
object concentrated in degree 0.
\end{Proof}

\begin{Def}{\rm
The {\it center} of a category $\AAA$ is the set of natural
transformations from $\Id_{\AAA}$ to itself. We denote the center of
a category $\AAA$ by $Z(\AAA)$. By $Z(\AAA)^*$ we denote the subset
of $Z(\AAA)$ formed by natural isomorphisms. If
$\theta:\AAA\rightarrow \BB$ is a functor, then $\alpha\in Z(\AAA)$
determines a natural transformation
$\theta(\alpha):\theta\rightarrow\theta$ by the formula
$\theta(\alpha)_x=\theta(\alpha_x)$. It is clear that if $\theta$ is
an equivalence, then any natural isomorphism from $\theta$ to
$\theta$ is of the form $\theta(\alpha)$ ($\alpha\in Z(\AAA)^*$). }
\end{Def}

Now let $F_{\theta}$ be an element of $\TP(\AAA,\BB)$. We want to
determine when $F_{\theta}$ lies in the image of $\Phi_{\AAA,\BB}$.
Let the group $G$ be given by generators and relations
$G=<\{a\}_{a\in A}|\{b\}_{b\in B}>$. We know from Proposition
\ref{sur_in} that $F_{\theta}\in\Im(\Phi_{\AAA,\BB})$ iff there is a
$G$-equivariance adjuster for $\theta$. In particular, if
$F_{\theta}\in\Im(\Phi_{\AAA,\BB})$, then there is some natural
isomorphism $\varphi_a:\theta a\rightarrow a\theta$ for any $a\in
A$.

Define $\varphi_{a^{-1}}:\theta a^{-1}\rightarrow a^{-1}\theta$ by
the formula
$$
\varphi_{a^{-1},y}={}^{a^{-1}}(\varphi_{a,a^{-1}y}^{-1}).
$$
Denote $\tilde A:=A\cup\{a^{-1}|a\in A\}$. Let us define natural
isomorphisms $\varphi_{a_1,\dots,a_n}:\theta a_1\dots a_n\rightarrow
a_1\dots a_n\theta$ for all families $a_1,\dots,a_n\in \tilde A$. We
have done this for the case $n=1$. Let
$\varphi_{a_1,\dots,a_{n-1}}:\theta a_1\dots a_{n-1}\rightarrow
a_1\dots a_{n-1}\theta$ be defined. Then we define
$\varphi_{a_1,\dots,a_n}$ by the formula
$$
\varphi_{a_1,\dots,a_n,x}={}^{a_1\dots a_{n-1}}\varphi_{a_n,x}\varphi_{a_1,\dots,a_{n-1},a_nx}.
$$
Let $a_1,\dots,a_n\in\tilde A$ be such elements that $a_1\dots
a_n\in B$. Then $\varphi_{a_1,\dots,a_n}$ is a natural isomorphism
from $\theta$ to itself. So there is a family
$\alpha=(\alpha_b)_{b\in B}$ of elements of $Z(\BB)^*$ such that
$\varphi_{a_1,\dots,a_n}=\theta(\alpha_b)$ for $b=a_1\dots a_n\in
B$. As it was mentioned above any $\sigma\in\Aut\BB$ induces an
automorphism $\sigma\in\Aut(\Kb\BB)$. It is clear that $\sigma$ is a
standard derived equivalence lying in $\Pic(\BB)$. Let
$\epsilon_a:\sigma a\rightarrow a\sigma$ ($a\in A$) be a family of
natural isomorphisms. We define $\epsilon_{a_1,\dots,a_n}:\sigma
a_1\dots a_n\rightarrow a_1\dots a_n\sigma$ for
$a_1,\dots,a_n\in\tilde A$ analogously to the definition of
$\varphi_{a_1,\dots,a_n}$. For $b=a_1\dots a_n\in B$ we define
$\epsilon_b$ by the formula $\epsilon_b=\epsilon_{a_1,\dots,a_n}$.

\begin{Def} {\rm In the above notation the family of isomorphisms
$\varphi=(\varphi_a)_{a\in A}$ is called an {\it approximate
equivariance adjuster} for $\theta$. The family $\alpha$ is called
an {\it equivariance error} for $\varphi$. The family
$\epsilon=(\epsilon_a)_{a\in A}$ is called an {\it equivariance
$\sigma$-correction} for $\alpha$ if $\alpha_{b,\sigma
y}\epsilon_{b,y}=\Id_{\sigma y}$ for any $b\in B$ and $y\in\BB$. }
\end{Def}

\begin{theorem}\label{PsiIm}
Suppose that $G$ is given by generators and relations. Let $\theta:\BB\rightarrow\Kb\AAA$ be a tilting functor. Suppose that $\varphi$ is an approximate equivariance
adjuster for $\theta$ with equivariance error $\alpha$. Let $\sigma$
be some automorphism of $\BB$. Then $F_{\theta}\sigma$ lies in the
image of $\Phi_{\AAA,\BB}$ iff there exists an equivariance
$\sigma$-correction for $\alpha$.
\end{theorem}
\begin{Proof}By Proposition \ref{sur_in} $F_{\theta}\sigma\in\Im\Phi_{\AAA,\BB}$
iff there is a $G$-equivariance adjuster for $\theta\sigma$. Suppose
that $\psi$ is a $G$-equivariance adjuster for $\theta\sigma$. Then
direct calculations show that $\epsilon$ defined by the equalities
$\psi_{a,x}=\varphi_{a,\sigma x}\circ\theta \epsilon_{a,x}$ is an
equivariance $\sigma$-correction for $\alpha$.

Assume now that $\epsilon$ is an equivariance $\sigma$-correction
for $\alpha$. For $a_1,\dots,a_n\in \tilde A$ we define
$\psi_{a_1\dots a_n}:\theta \sigma a_1\dots a_n\rightarrow a_1\dots
a_n\theta\sigma$ by the formula
$$
\psi_{a_1\dots a_n,x}=\varphi_{a_1,\dots,a_n,\sigma x}\circ\theta\epsilon_{a_1,\dots,a_n,x}.
$$
The correctness of this definition follows from the fact that
$\epsilon$ is an equivariance $\sigma$-correction for $\alpha$. It
can be verified by direct calculations that $\psi$ is a
$G$-equivariance adjuster for $\theta\sigma$.
\end{Proof}

\section{$G$-grading and the image of $\Psi$}

In this section we give a description of the image of
$\Psi_{\AAA,\BB}$ in the case where $G$ is a finite group. We need
the finiteness of $G$ to prove the following lemma.

\begin{Lemma}\label{tiltcrit}
Let $G$ be a finite group and $\XX$ be a subcategory of $\Kb\AAA$
such that for any $U\in\XX$ and any $g\in G$ there is $V\in\XX$ such
that ${}^gU\cong V$. Then $\XX$ is a tilting subcategory for $\AAA$
iff $P_{\bullet}\XX$ is a tilting subcategory for $\AAA/G$.
\end{Lemma}
\begin{Proof}
For the proof of the fact that $P_{\bullet}\XX$ is a tilting
subcategory for $\AAA/G$ if $\XX$ is a tilting subcategory for
$\AAA$  see the proof of \cite[Theorem 4.7]{Asashiba}.

Now let $P_{\bullet}\XX$ be a tilting subcategory for $\AAA/G$. Let
us prove that $\Kb\AAA(U,V[i])=0$ for $U,V\in\XX$, $i\not=0$.
Suppose that it is not true. Then there are $U,V\in\XX$ and
$i\not=0$ such that $\Kb\AAA(U,V[i])\not=0$. Let $f$ be a nonzero
element of $\Kb\AAA(U,V[i])$. Then $P_{\bullet}f$ is  a nonzero
element of $\Kb\AAA/G(P_{\bullet}U,P_{\bullet}V[i])$ and so
$P_{\bullet}\XX$ is not a tilting subcategory for $\AAA/G$. It
remains to prove that any representable functor lies in the
subcategory $\thick\XX$. Note that if $|G|<\infty$, then
$P^{\bullet}$ sends finitely generated modules to finitely generated
modules and so it induces a functor
$P^{\bullet}:\Kb(\AAA/G)\rightarrow\Kb\AAA$. Let us consider a
functor $\AAA(-,x)$ ($x\in\AAA$). By our assumption $\AAA(-,Px)\cong
P_{\bullet}\big(\AAA(-,x)\big)$ lies in $\thick P_{\bullet}\XX$. Let
us consider the subcategory $\thick P^{\bullet}P_{\bullet}\XX$ of
$\Kb\AAA$. It contains the subcategory $P^{\bullet}\thick
P_{\bullet}\XX$ and so contains
$P^{\bullet}P_{\bullet}\big(\AAA(-,x)\big)\cong\oplus_{g\in
G}\AAA(-,gx)$. Since $\thick P^{\bullet}P_{\bullet}\XX$ is closed
under direct summands it contains $\AAA(-,x)$. It remains to prove
that $\thick\XX$ contains $P^{\bullet}P_{\bullet}U$ for any
$U\in\XX$. But $\thick\XX$ contains ${}^gU$ for any $U\in\XX$ and
any $g\in G$. So it contains
$P^{\bullet}P_{\bullet}U\cong\oplus_{g\in G}{}^gU$ for any
$U\in\XX$.
\end{Proof}

\begin{Def}{\rm
A {\it $G$-graded category} is a category $\AAA$ having a family of
direct sum decompositions $\AAA(x,y)=\oplus_{g\in G}\AAA(x,
y)^{(g)}$ $(x, y \in \AAA)$ of $\kk$-modules such that the
composition of morphisms gives the inclusions $\AAA(y,
z)^{(g)}\AAA(x, y)^{(h)}\subset\AAA(x, z)^{(gh)}$ for all $x, y, z
\in\AAA$ and $g,h\in G $. If $f \in \AAA(x, y)^{(g)}$, then we say
that $f$ is of $G$-degree $g$. A functor $F: \AAA \rightarrow \AAA'$
between $G$-graded categories is called {\it degree-preserving} if
$F(\AAA(x, y)^{(g)}) \subset \AAA'(Fx, Fy)^{(g)}$ for all $x, y \in
\AAA$ and $g\in G$. }
\end{Def}

By \cite[Lemma 5.4]{Asashiba} there is a $G$-grading on $\BB/G$ such
that $Pf$ is of degree $1_G$ and $s_{g,x}$ is of degree $g^{-1}$.

\begin{Def}{\rm
Let $\AAA$ be a $G$-graded category. A {\it $G$-graded
$\AAA$-complex} is an $\AAA$-complex $U$ with a family of direct sum
decompositions $U(x)=\oplus_{g\in G}U(x)^{(g)}$ ($x\in\AAA$) such
that $d_{U,x}(U(x)^{(g)})\subset U(x)^{(g)}$ and
$U(f)(U(x)^{(g)})\subset U(y)^{(gh)}$ for all $x,y\in\AAA$, $g,h\in
G$, $f\in\AAA(y,x)^{(h)}$. If $U,V$ are $G$-graded complexes, then
we say that $f:U\rightarrow V$ is of degree $h$ if
$f_x(U(x)^{(g)})\subset V(x)^{(hg)}$. Let us now define a category
$\KbG\AAA$. It's objects are $G$-graded $\AAA$-complexes $U$ which
lie in $\Kb\AAA$ considered as complexes without grading. If
$U,V\in\KbG\AAA$, then $\KbG\AAA(U,V)=\Kb\AAA(U,V)$. }
\end{Def}

Let $\AAA$ be a $G$-graded category and $U,V\in\KbG\AAA$. Then we
denote by $\KbG\AAA(U,V)^{(g)}$ the set of morphisms in
$\Kb\AAA(U,V)$ which can be presented by a morphism of degree $g$ in
$\CC\AAA(U,V)$. It is not hard to check that
$\KbG\AAA(x,y)=\oplus_{g\in G}\KbG\AAA(U, V)^{(g)}$ and this
decomposition turns $\KbG\AAA$ into a $G$-graded category. For
$U\in\KbG\AAA$ we denote by $\bar U$ the corresponding object of
$\Kb\AAA$. Note that any equivalence
$\theta:\BB\rightarrow \KbG\AAA$ determines an equivalence
$\bar\theta:\BB\rightarrow \Kb\AAA$ in an obvious way.

\begin{theorem}
Suppose that $G$ is a finite group. Let
$F:\Kb(\BB/G)\rightarrow\Kb(\AAA/G)$ be a standard equivalence. Then
$F$ lies in the image of $\Psi_{\AAA,\BB}$ iff there is a
degree-preserving functor $\mu:\BB/G\rightarrow \KbG(\AAA/G)$ such
that $\bar\mu$ is a tilting functor and $F\cong F_{\bar\mu}$.
\end{theorem}

\begin{Proof}
If $F$ lies in the image of $\Psi_{\AAA,\BB}$, then it is isomorphic
to $F_{\mu_{\theta,\psi}}$ for some tilting functor
$\theta:\BB\rightarrow\Kb\AAA$ and a $G$-equivariance adjuster
$\psi$ for $\theta$. Let us define a $G$-grading on
$\mu_{\theta,\psi}(Py)=P_{\bullet}\theta(y)$ ($y\in\BB$) as follows:
$(P_{\bullet}\theta(y))(Px)^{(g)}=\theta(y)(gx)$. Then it can be
easily verified that $\mu_{\theta,\psi}$ defines a degree-preserving
functor from $\BB/G$ to $\KbG(\AAA/G)$. Note that this part of the
prove does not require the finiteness of $G$.

Now let $\mu:\BB/G\rightarrow\KbG(\AAA/G)$ be a degree-preserving
functor such that $\bar\mu$ is a tilting functor. Let us prove that
$F_{\bar\mu}$ lies in the image of $\Psi_{\AAA,\BB}$. By Proposition
\ref{sur_in} it is enough to find a tilting functor
$\theta:\BB\rightarrow \Kb\AAA$ and a $G$-equivariance adjuster
$\psi$ for $\theta$ such that $F_{\bar\mu}\cong
F_{\mu_{\theta,\psi}}$.

For $U\in\KbG(\AAA/G)$ we define $\widetilde U\in\Kb\AAA$ in the
following way. It is defined on objects by the formula $\widetilde
U(x)=U(Px)^{(1)}$ and on morphisms by the formula $\widetilde
U(f)=U(Pf)|_{U(x)^{(1)}}$ ($f\in\AAA(x,y)$). The differential
$d_{\widetilde U}$ is defined by the formula $d_{\widetilde
U,x}=d_{U,Px}|_{U(Px)^{(1)}}$. The correctness of this definition
follows from the definition of a $G$-graded complex. Also we define
a morphism $\xi_U:\bar U\rightarrow P_{\bullet}\widetilde U$ as
follows:
$$
\xi_{U,x}=\bigoplus_{g\in G}U(s_{g,x})|_{U(Px)^{(g)}}:\bigoplus_{g\in G}U(Px)^{(g)}\rightarrow \bigoplus_{g\in G}U(Pgx)^{(1)}.
$$

Let us now define $\theta:\BB\rightarrow \Kb\AAA$. We define it on
objects by $\theta(y)=\widetilde{\mu(Py)}$. For $f\in\BB(y,z)$ we
define the natural transformation $\theta(f)$ by the formula
$$
\theta(f)_x=\mu(Pf)_{Px}|_{\mu(Py)(Px)^{(1)}}:\theta(y)(x)\rightarrow\theta(z)(x)
$$
for all $x\in\AAA$. Let $\psi_{g,y,x}:\theta(gy)(x)\rightarrow {}^g\theta(y)(x)$ ($x\in\AAA$, $y\in\BB$, $g\in G$) be the composition
$$
\theta(gy)(x)=\mu(Pgy)(Px)^{(1)}\stackrel{\mu(s_{g,y})_x}\longrightarrow \mu(Py)(Px)^{(g^{-1})}\stackrel{\mu(y)(s_{g^{-1},x})}\longrightarrow\mu(Py)(Pg^{-1}x)^{(1)}={}^g\theta(y)(x).
$$

It is not hard to prove that the following conditions hold:\\
1) $\xi_U$ is an isomorphism in $\Kb(\AAA/G)$;\\
2) $\theta$ induces an equivalence from $\BB$ to $\theta(\BB)$;\\
3) $\psi$ is a $G$-equivariance adjuster for $\theta$;\\
4) The family of morphisms $\zeta_{Py}=\xi_{\mu(Py)}:\overline{\mu(Py)}\rightarrow P_{\bullet}\theta(y)=\mu_{\theta,\psi}(Py)$ defines an isomorphism from $\bar\mu$ to $\mu_{\theta,\psi}$.\\
It follows from 1)--3) and Lemma \ref{tiltcrit} that $\theta(\BB)$
is a tilting subcategory for $\AAA$, hence the theorem is proved.
\end{Proof}

\section{The action of the Nakayama automorphism on Frobenius algebras}

From here on we assume that $\kk$ is a field. Let $R$ be an
associative finite dimensional $\kk$-algebra. It can be considered
as a category with one object and so the results of the previous
sections can be applied to $R$. We denote the unique object of $R$
by $e$. We apply these results in the case where $R$ is a Frobenius
algebra and $G$ is a finite cyclic group which acts on $R$ by powers
of a Nakayama automorphism. First, let us recall the definition of a
Frobenius algebra.

\begin{Def}{\rm
An algebra $R$ is called {\it Frobenius} if there is a linear map
$\epsilon:R\rightarrow k$ such that the bilinear form $\langle
a,b\rangle=\epsilon(ab)$ is nondegenerate. The Nakayama automorphism
$\nu:R\rightarrow R$ is the automorphism which satisfies the
equation $\langle a,b\rangle=\langle b,\nu(a)\rangle$ for all
$a,b\in R$. If the bilinear form on $R$ can be chosen in such a way
that $\langle a,b\rangle=\langle b,a\rangle$ for all $a,b\in R$,
then the algebra $R$ is called {\it symmetric}. }
\end{Def}

From here on we fix some Frobenius algebra $R$, and some Nakayama
automorphism $\nu$ of $R$. Moreover, we assume that there is an
integer $n>0$ such that $\nu^n=\Id_R$. Then the cyclic group
$G=<g|g^n>$ acts on $R$ by the following rule: $\Delta(g)=\nu$. Note
that $R/G$ is a symmetric algebra. Indeed, if $\langle \mbox{ },
\mbox{ }\rangle$ is the bilinear form on $R$, then $\langle \mbox{
}, \mbox{ }\rangle_G$ defined by the formula
$$
<s_{g^k}Pa,s_{g^l}Pb>_G=\delta_{g^{k+l+1},1_G}<\nu^la,b>
$$
is the desired bilinear form on $R/G$. So the maps $\Phi_R$ and
$\Psi_R$ allow us to transfer some information from the derived
Picard group of a symmetric algebra to the derived Picard group of a
Frobenius algebra.

We have the following application of Theorem \ref{PsiIm}.

\begin{prop}\label{Psisur1}
$1)$ If for any $a\in Z(R)^*$ there exists an element $b\in Z(R)^*$ such that $ab^n=1$, then $\Phi_R$ is surjective.\\
$2)$ If for any $a\in Z(R)^*$ there exists an element $b\in R^*$ and
an automorphism $\sigma\in\Aut R$ such that
$b\nu(b)\dots\nu^{n-1}(b)=a$ and
$\sigma\nu\sigma^{-1}(c)=b\nu(c)b^{-1}$ for any $c\in R$, then
$\Cok\Phi_R$ is generated by the images of elements from $\Pic(R)$.
\end{prop}
\begin{Proof}
Note that the functor ${}^{\nu}(-):\Kb R\rightarrow \Kb R$ is
isomorphic to the Nakayama functor. Then by \cite[Proposition
5.2]{Rickard} we have an isomorphism $\eta_F:F\circ\nu\cong\nu\circ
F$ for any standard equivalence $F:\Kb R\rightarrow \Kb R$. If $F$
is given by a tilting functor $\theta:R\rightarrow \Kb R$, then
$\theta\cong F\YY_R$. So the isomorphism $\eta_{F,\YY(e)}\circ
F(\phi_{g,e}):F\YY_R\nu(e)\rightarrow \nu F\YY_R(e)$ gives an
isomorphism $\varphi_{g,e}:\theta \nu(e)\rightarrow \nu\theta(e)$.
Then $\varphi$ is an approximate equivariance adjuster for $\theta$.

1) Note that $\nu(b)=b$ for any $b\in Z(R)$. Then it follows from
the assumption that there is an equivariance $\Id_R$-correction for
any equivariance error. So $\Phi_R$ is surjective by Theorem
\ref{PsiIm}.

2) It follows from the assumption that for any equivariance error
$a$ there is some $\sigma\in\Aut R$ such that there exists an
equivariance $\sigma$-correction for $a$. So by Theorem \ref{PsiIm}
for any $F\in \TP(R)$ there is some $\sigma\in\Aut R$ such that
$F\sigma$ lies in the image of $\Phi_R$. Then the image of $F$ in
$\Cok\Phi_R$ equals the image of $\sigma^{-1}$. So the assertion
follows from the fact that $\sigma^{-1}\in\Pic(R)$.

\end{Proof}

\begin{coro}\label{sur_coro}
If the field $\kk$ is algebraically closed and it's characteristic
does not divide $n$, then $\Phi_R$ is surjective.
\end{coro}
\begin{Proof}
We may assume that $R$ is an indecomposable algebra. Then any
element $a\in Z(R)^*$ is of the form $a=\kappa(1+Q)$ for some
$\kappa\in\kk$ and nilpotent $Q\in Z(R)$. Since $\kk$ is
algebraically closed, there is some $\bar\kappa\in\kk$ such that
$\kappa=\bar\kappa^n$.  Then $a=b^n$ for
$$
b=\bar\kappa\sum\limits_{i\ge 0}\frac{\prod\limits_{j=0}^{i-1}\left(\frac{1}{n}-j\right)}{i!}Q^i.
$$
\end{Proof}

\section{Application: generators of the derived Picard group of a self-injective Nakayama algebra}

From here on we assume that $\kk$ is an algebraically closed field.
In this section we apply the methods of  the previous sections to
obtain generators of the derived Picard group of algebras $\NN(nm,
tm)$ defined in the following way. Let $m,n,t>0$ be some integers.
We suppose that $n$ and $t$ are coprime. Let $\QQ(nm)$ be a cyclic
quiver with $nm$ vertices, i.e. the quiver whose vertex set is
$\mathbb{Z}_{nm}$ and whose arrows are $\beta_i:i\rightarrow i+1$
($i\in\mathbb{Z}_{nm}$). Let $\II(nm,tm)$ be an ideal in the path
algebra of $\QQ(nm)$ generated by all paths of length $tm+1$. We
denote $\NN(nm, tm):=\kk\QQ(nm)/\II(nm,tm)$. For
$i\in\mathbb{Z}_{nm}$ we denote by $e_i$ the primitive idempotent
corresponding to the vertex $i$ and by $P_i$ the projective module
$e_i\NN(nm, tm)$. For a path $w$ from the vertex $i$ to the vertex
$j$ we denote by $w$ the unique homomorphism from $P_i$ to $P_j$
which sends $e_i$ to $w$ as well. Also we introduce the following
auxiliary notation:
$$
\beta_{i,k}=\beta_{i+k-1}\dots\beta_i.
$$
It is well-known that $\NN(nm, tm)$ is a Frobenius algebra with a
Nakayama automorphism $\nu$ defined as follows: $\nu(e_i)=e_{i-tm}$
and $\nu(\beta_i)=\beta_{i-tm}$. If $U$ is a module, then we also
denote by $U$ the corresponding complex concentrated in degree 0.
For $i\in\mathbb{Z}_{nm}$, $1\le k\le m-1$ we introduce the complex
$$
X_i:=P_{i-tm}\stackrel{\beta_{i-tm}}\rightarrow P_{i-tm+1}\stackrel{\beta_{i-tm+1,tm}}\rightarrow P_{i+1}
$$
concentrated in degrees -2, -1 and 0 and the complexes
$$
Y_{i,k}:=P_i\stackrel{\beta_{i,k}}\rightarrow P_{i+k}
$$
concentrated in degrees 0 and 1.

If $m>1$, then for $0\le l\le m-1$ we introduce the $\NN(nm, tm)$-complex
$$
H_l^{nm}=\Big(\bigoplus_{i\in\mathbb{Z}_{nm}, m\nmid i-l}P_i\Big)\oplus\Big(\bigoplus_{i\in\mathbb{Z}_{nm}, m\mid i-l}X_i\Big).
$$
In this case we can define an algebra isomorphism
$$
\theta_l^{nm}:\NN(nm, tm)\rightarrow \Kb(\NN(nm,
tm))(H_l^{nm},H_l^{nm}).
$$
We define it on idempotents by the formula
$$
\theta_l^{nm}(e_i)=\begin{cases}
\Id_{P_i}&\mbox{ if $m\nmid i-l$ and $m\nmid i-1-l$,}\\
\Id_{P_{i+1}}&\mbox{ if $m\mid i-l$,}\\
\Id_{X_i}&\mbox{ if $m\mid i-1-l$.}
\end{cases}
$$
We define $\theta_l^{nm}(\beta_i)$ ($i\in\mathbb{Z}_{nm}$) in the
obvious way (it equals 0 in all degrees except for the zero degree
and equals $\beta_i$, $\beta_i\beta_{i+1}$ or $\Id_{P_{i+1}}$
depending on $i$ in the zero degree).

If $m>1$ and $t=1$, then for $0\le l\le m-1$ we introduce the $\NN(nm, m)$-complex
$$
Q_l^{nm}=\bigoplus_{i\in\mathbb{Z}_{nm}, m\mid i-l}\Big(P_i\oplus \bigoplus_{k=1}^{m-1}Y_{i,k}\Big).
$$
In this case we can define an algebra isomorphism
$$
\varepsilon_l^{nm}:\NN(nm, m)\rightarrow \Kb(\NN(nm, m))(Q_l^{nm},Q_l^{nm}).
$$
We define it on idempotents by the fomula
$$
\varepsilon_l^{nm}(e_i)=\begin{cases}
\Id_{P_i}&\mbox{ if $m\mid i-l$,}\\
\Id_{Y_{i+k,m-k}}&\mbox{ if $m\mid i+k-l$ for some $k$, $1\le k\le
m-1$.}
\end{cases}
$$
We define $\varepsilon_l^{nm}(\beta_i)$ ($i\in\mathbb{Z}_{nm}$) in the following way. It equals 0 in all degrees except 0 and 1. In degree 0 it equals $\Id_{P_{i+k}}$ if  $m\mid i+k-l$ for $1\le k\le m-1$ and equals $\beta_{i,m}$ if $m\mid i-l$. In degree 1 it equals $\beta_{i+m}$ if  $m\nmid i-l$ and $m\nmid i+1-l$ and equals 0 if $m\mid i-l$ or $m\mid i+1-l$.

In this section we prove the following theorem.

\begin{theorem}\label{Nakayama_gen}
$1)$ If $m=1$, then $\TP(\NN(n, t))$ is generated by the shift and $\Pic(\NN(n, t))$;\\
$2)$ If $m>1$, $t>1$, then $\TP(\NN(nm, tm))$ is generated by the shift, $\Pic(\NN(nm, tm))$ and $F_{\theta_l^{nm}}$ ($l\in\mathbb{Z}_m$);\\
$3)$ If $m>1$, $t=1$, then $\TP(\NN(nm, m))$ is generated by the shift, $\Pic(\NN(nm, m))$, $F_{\theta_l^{nm}}$ and $F_{\varepsilon_l^{nm}}$ ($l\in\mathbb{Z}_m$).
\end{theorem}

It is clear that $\NN(nm, tm)$ is symmetric iff $n=1$. If in
addition $m=1$, then $\NN(1,t)$ is a local algebra and so
$\TP(\NN(1, t))$ is generated by the shift and $\Pic(\NN(1, t))$ by
the results of \cite{Yek}, \cite{bla-bla-bla}. So the first
assertion of Theorem \ref{Nakayama_gen} holds for $n=1$. The
assertions 2) and 3) of the theorem for $n=1$ follow from the
results of \cite{Zvonareva}, where the set of generators of the
derived Picard group was described in the case $n=1$, $m>1$.
Moreover, it was proved there that any element of $\TP(\NN(m, tm))$
is of the form $UV$, where $V\in\Pic(\NN(m, tm))$ and $U$ is a
product of elements listed in the points 2)--3) except for
$\Pic(\NN(m, tm))$.

Now let us consider $n>1$. Let $G=\langle g\mid g^n\rangle$ be a
cyclic group which acts on $\NN(nm, tm)$ by the rule
$\Delta(g)=\nu$. It is well-known that $\NN(nm, tm)/G$ is Morita
equivalent to $\NN(m, tm)$. We need the explicit formula for this
equivalence to obtain the isomorphism of the derived Picard groups
defined by it. Let $W=\oplus_{j\in\mathbb{Z}_n}W_j$ where $W_j$ is
isomorphic to $\NN(m, tm)$ as a right $\NN(m, tm)$-module. As it was
mentioned above every path $w$ in $\QQ(m)$ defines a homomorphism
$w:W_j\rightarrow W_j$. Thus, there is a left $\NN(m, tm)$-module
structure on $W_j$. Let us define a left $\NN(nm, tm)/G$-module
structure on $W$. Let $s_{i,j}:W_i\rightarrow W_j$
($i,j\in\mathbb{Z}_n$) be the isomorphism arising from $\Id_{\NN(m,
tm)}$. Let $i\in\mathbb{Z}_{nm}$ be represented by an integer number
$0\le \bar i\le nm-1$. Present $\bar i$ in the form $\bar i=\bar
qm+\bar r$, where $0\le\bar r< m$. Let $q\in\mathbb{Z}_n$ and
$r\in\mathbb{Z}_m$ be elements represented by $\bar q$ and $\bar r$
respectively. Consider an element $x\in W$. Suppose that $x\in W_j$
for some $j\in\mathbb{Z}_n$. Then we define
$$
(Pe_i)x=(\delta_{q,j}e_{r})x, (P\beta_i)x=(\delta_{q,j}\beta_{r})x\mbox{ and }s_{g^l}x=s_{j,j+ltm}(x).
$$
It is clear that in such a way $W$ becomes a $\NN(nm, tm)/G-\NN(m,
tm)$-bimodule which induces a Morita equivalence
$$T_W=-\ot_{\NN(nm, tm)/G} W:\Kb(\NN(nm, tm)/G)\rightarrow\Kb(\NN(m,tm)).$$
 We define $L:\TP(\NN(m,tm))\rightarrow \TP(\NN(nm, tm)/G)$ by the formula
$$L(F)=\bar T_W\circ F\circ T_W,$$
where $\bar T_W$ is a quasi-inverse equivalence for $T_W$.
 It is clear that $L$ sends the shift to the shift and $\Pic(\NN(m,tm))$ to $\Pic(\NN(nm, tm)/G)$.

There are $G$-equivariance adjusters $\psi_l^{nm}$ for $\theta_l^{nm}$ and $\varphi_l^{nm}$ for $\varepsilon_l^{nm}$. The maps $\psi_{l,g^p}^{nm}$ and $\varphi_{l,g^p}^{nm}$ can be constructed in the obvious way as the sums of the isomorphisms of the form
$$
P_{i-ptm}\cong {}^{\nu^p}P_i,X_{i-ptm}\cong {}^{\nu^p}X_i\mbox{ and }Y_{i-ptm,k}\cong {}^{\nu^p}Y_{i,k}.
$$
Then we have the following lemma.

\begin{Lemma}\label{move_gen}
$1)$ If $m>1$, then for $0\le l\le m-1$ we have
$
L(F_{\theta_l^m})\cong F_{\mu_{\theta_l^{nm},\psi_l^{nm}}}.
$

$2)$If $m>1$ and $t=1$, then for $0\le l\le m-1$ we have
$
L(F_{\varepsilon_l^m})\cong F_{\mu_{\varepsilon_l^{nm},\varphi_l^{nm}}}.
$
\end{Lemma}
\begin{Proof}
1) We will prove that
\begin{equation}\label{eq_0}
F_{\theta_l^m}\circ T_W\cong T_W\circ F_{\mu_{\theta_l^{nm},\psi_l^{nm}}}.
\end{equation}
Let us describe the left part of this equality. Let
$H=\bigoplus_{j\in\mathbb{Z}_n}H_j$, where $H_j\cong H_l^m$ as a
right $\NN(m,tm)$-complex. Denote by $s'_{i,j}:H_i\rightarrow H_j$
($i,j\in\mathbb{Z}_n$) the isomorphism arising from $\Id_{H_l^m}$.
In addition, for $u\in\Kb(\NN(m,tm))(H_l^m,H_l^m)$ we denote by $u$
the corresponding morphism from $H_j$ to $H_j$. Let us define
$\theta:\NN(nm, tm)/G\rightarrow \Kb(\NN(m,tm))(H,H)$. Consider
$i\in\mathbb{Z}_{nm}$. Let $q\in\mathbb{Z}_n$ and $r\in
\mathbb{Z}_m$ be as above. Then for $x\in H_j$ ($j\in\mathbb{Z}_n$)
we define
 $$
 \theta(e_i)(x)=\delta_{j,q}\theta_l^m(e_r)(x), \theta(\beta_i)(x)=\delta_{j,q}\theta_l^m(\beta_r)(x)\mbox{ and } \theta(s_{g^l})(x)=s'_{j,j+ltm}(x).
 $$
Then the left part of the equality \eqref{eq_0} gives $F_{\theta}$.
It is not hard to construct an isomorphism
 $$\xi:H\rightarrow (P_{\bullet}H_l^{nm})\ot_{\NN(nm, tm)/G} W$$
such that
$\xi\theta(c)=(\mu_{\theta_l^{nm},\psi_l^{nm}}(c)\ot_{\NN(nm,
tm)/G}\Id_W)\xi$ for any $c\in\NN(nm, tm)/G$. The existence of such
$\xi$ gives the isomorphism \eqref{eq_0}.

2) The proof is similar and so it is left to the reader.
\end{Proof}

Let us now apply the results of the previous sections to the algebra
$\NN(nm,tm)$.

\begin{Lemma}\label{Nakayama_sur}
$1)$ If $\charr\kk\nmid n$, then $\Phi_{\NN(nm, tm)}$ is surjective.\\
$2)$ If $\charr\kk\mid n$, then $\Cok\Phi_{\NN(nm, tm)}$ is generated by images of elements from $\Pic(\NN(nm, tm))$.
\end{Lemma}
\begin{Proof}
1) Follows directly from Corollary \ref{sur_coro}.

2) Note, that in this case $n>1$. It is enough to prove that the
condition of the second part of Proposition \ref{Psisur1} is
satisfied. Consider $a\in Z(\NN(nm, tm))^*$. Let us introduce the
notation $u:=\sum_{i\in\mathbb{Z}_{nm}}\beta_{i,nm}$, $\bar
u:=\sum_{i=1}^{m}\beta_{i,nm}$. It can be easily proved that
$a=\sum_{0\le k\le \frac{t}{n}}c_ku^k$ for some $c_k\in\kk$,
$c_0\not=0$. Since $\kk$ is algebraically closed we may assume that
$c_0=1$. Then $a=b\nu(b)\dots\nu^{n-1}(b)$ for $b=1+\sum_{1\le k\le
\frac{t}{n}}c_k\bar u^k$. Let us denote by $\gamma$ the automorphism
of $\NN(nm, tm)$ defined by the formula $\gamma(x)=b\nu(x)b^{-1}$
for $x\in \NN(nm, tm)$. It remains to find such
$\sigma\in\Aut\NN(nm, tm)$ that $\sigma^{-1}\nu\sigma=\gamma$.

The automorphism $\gamma$ equals $\nu$ on the idempotents and is
defined on the arrows by the formula
$$
\gamma(\beta_i)=\begin{cases}
\beta_{i-tm},&\mbox{ if $i\not=tm$ and $i\not=(t+1)m$},\\
a\beta_0,&\mbox{ if $i=tm$},\\
a^{-1}\beta_{m},&\mbox{ if $i=(t+1)m$}.
\end{cases}
$$
Let $0<p<n$ be such a number that $n\mid pt-1$. We define $\sigma$
in the following way. It is identical on the idempotents and is
defined on arrows in such a way that
$$
\sigma^{-1}(\beta_i)=\begin{cases}
a^{-1}\beta_i,&\mbox{ if $i=(1+k)tm$ for some $0\le k< p$},\\
\beta_i,&\mbox{ otherwise}.
\end{cases}
$$
It is easy to verify that $\sigma^{-1}\nu\sigma=\gamma$.
\end{Proof}

\begin{Proof}[Proof of Theorem \ref{Nakayama_gen}] It was mentioned above that the theorem is true for $n=1$.
Consider $n>1$. We want to prove that some set of elements of
$\TP(\NN(nm, tm))$ generates it. Denote this set by $M$. It follows
from Lemma \ref{Nakayama_sur} that the image of $\Phi_{\NN(nm, tm)}$
and $\Pic(\NN(nm, tm))$ generates $\TP(\NN(nm, tm))$.

Let $F\in\TP_G(\NN(nm, tm))$. It is enough to prove that $\Phi_{\NN(nm, tm)}(F)$ lies in the subgroup of $\TP(\NN(nm, tm))$ generated by $M$.
It follows from Lemma \ref{move_gen} and the arguments above that $\Psi_{\NN(nm, tm)}(F)=\Psi_{\NN(nm, tm)}(U)V$ for some $V\in\Pic(\NN(m, tm))$ and some $U\in\TP_G(\NN(nm, tm))$ such that $\Phi_{\NN(nm, tm)}(U)$ lies in the subgroup generated by $M$. By Theorem \ref{Picard} we have $\Phi_{\NN(nm, tm)}(U^{-1}F)\in\Pic(\NN(nm, tm))$. Then $\Phi_{\NN(nm, tm)}(U^{-1}F)$ lies in the subgroup generated by $M$ and, consequently, $\Phi_{\NN(nm, tm)}(F)$ lies in the subgroup generated by $M$.
\end{Proof}


\begin{thebibliography}{99}
\bibitem{Keller}
B. Keller, {\it Deriving DG categories}. ---  Ann. Sci. \'Escole Norm. Sup. 27, No. 4, 63--102 (1994).

\bibitem{Asashiba}
H. Asashiba, {\it A generalization of Gabriel's Galois covering functors and derived equivalences}. ---  J. Algebra 334, No. 1, 109--149 (2011).

\bibitem{Rickard}
J. Rickard, {\it Derived equivalences as derived functors}. ---  J. London Math. Soc. 43, 37--48 (1991).

\bibitem{Zvonareva}
A. Zvonareva, {\it On the derived Picard group of the Brauer star algebra}. ---  arXiv:1401.6952.

\bibitem{LeChen}
J. Le, X.-W. Chen, {\it Karoubianness of a triangulated category}. --- J. Algebra 310, 452--457 (2007).

\bibitem{Rickard2}
J. Rickard, {\it Morita theory for derived categories}. ---  J. London Math. Soc. 39, 436--456 (1989).

\bibitem{Asashiba2}
H. Asashiba, {\it A generalization of Gabriel's Galois covering functors II: 2-categorical Cohen-Montgomery duality}. ---  arXiv:0905.3884.

\bibitem{Yek} A. Yekutieli, {\it Dualizing complexes, Morita equivalence and the
derived Picard group of a ring}. --- J. London Math. Soc. 60.3,
723--746 (1999).

\bibitem{bla-bla-bla}
R. Rouquier, A. Zimmermann, {\it Picard groups for derived module
categories}. --- Proc. London Math. Soc. 87.01, 197--225 (2003).
\end{thebibliography}
\end{document}